\documentclass{article}
\usepackage{amsmath}
\usepackage{amssymb}
\usepackage{amsthm}
\usepackage{float}
\newtheorem{Theorem}{Theorem}
\newtheorem{Conjecture}{Conjecture}
\newtheorem{Definition}{Definition}

\newtheorem{Lemma}{Lemma}

\usepackage{paralist}
\usepackage{graphicx}

\newcommand{\Includegraphics}[2]%
{\ifFigs\includegraphics[      width=#1]{#2}%
	\else  \includegraphics[draft,width=#1]{#2}\fi}

\begin{document}
	
	\title{Abundance of weird quasiperiodic attractors \\ 
	in piecewise linear discontinuous maps}
	\author{\textit{Laura Gardini}$^{1,2},$ 
		\textit{ Davide Radi}$^{2,3}$, 
		\textit{Noemi Schmitt}$^{4}$, \\ 
		\textit{Iryna Sushko}$^{3,5}$, 
		\textit{Frank Westerhoff}$^{4}$\\ 
		$^{1}${\small Dept of Economics, Society and Politics, University of Urbino
			Carlo Bo, Italy}\\ 
		$^{2}${\small Dept of Finance, V\v{S}B - Technical University of Ostrava,
			Ostrava, Czech Republic}\\ 
		$^{3}${\small Dept of Mathematics for Economic, Financial and Actuarial
			Sciences,} \\ {\small Catholic University of Milan, Italy} \\ 
		$^{4}${\small Dept of Economics, University of Bamberg, Germany}\\ 
		$^{5}${\small Inst. of Mathematics, NAS of Ukraine, Kyiv, Ukraine}}
	
	\date{}
	\maketitle
	
	\textbf{Abstract}
	\smallskip

In this work, we consider a class of $n$-dimensional, $n\geq2$, piecewise
linear discontinuous maps that can exhibit a new type of attractor, called a
\textit{weird quasiperiodic attractor}. While the dynamics associated with
these attractors may appear chaotic, we prove that chaotic attractors cannot occur.  
The considered class of $n$-dimensional maps allows for any finite number of
partitions, separated by various types of discontinuity sets. The key
characteristic, beyond discontinuity, is that all functions defining the map
have the same real fixed point. These maps cannot have hyperbolic cycles other
than the fixed point itself. We consider the two-dimensional case in detail.
We prove that in nongeneric cases, the restriction, or the first return, of
the map to a segment of straight line issuing from the fixed point 
is reducible to a piecewise linear
circle map. The generic attractor, different from the fixed point, is a weird
quasiperiodic attractor, which may coexist with other attractors or attracting
sets. We illustrate the existence of these attractors through numerous
examples, using functions with different types of Jacobian matrices, as well
as with different types of discontinuity sets. In some cases, we describe
possible mechanisms leading to the appearance of these attractors. We also
give examples in the three-dimensional space. Several properties of this new
type of attractor remain open for further investigation.$\medskip$

Keywords: Piecewise linear discontinuous maps; Weird quasiperiodic attractors; 
Piecewise linear circle maps; Attractors without hyperbolic cycles.

\section{Introduction}

The existence of a new type of attractor in two-dimensional (2D) discontinuous
piecewise linear (PWL) maps, which appear chaotic but are not, has been
observed in 2D PWL homogeneous systems that model economic dynamics
(see \cite{GRSSW-24,GRSSW-25c}). In a recent paper, \cite{GRSSW-25a},
we identified specific regions in the parameter space associated with this
type of attractor, called weird quasiperiodic attractor (WQA), focusing on a
2D discontinuous PWL homogeneous map derived from the well-known 2D border
collision normal form \cite{Simpson-14,Simpson-17}).

The goal of this work is to show that WQAs can be observed in a broad class of
2D PWL discontinuous maps, regardless of the type of borders of the
partitions where different functions are defined. We also examine particular
nongeneric cases where the attracting sets can be analyzed via a
one-dimensional (1D) restriction of the map or a first return map, ultimately
leading to a PWL circle map. These nongeneric cases have been recently
investigated in \cite{GRSSW-25b}, where we describe the dynamics of the
related class of 1D maps. Furthermore, we show that WQAs also occur in
three-dimensional maps and, more generally, may exist in $n$-dimensional
($n$D) maps.

To introduce a number of basic concepts, we begin with a 2D discontinuous PWL
map (often referred to as piecewise affine), in which the functions are
defined in two partitions and have the same real fixed point 
$(x,y)=(-\xi,-\eta).$ The system is given by:
\begin{equation}
M:\left\{
\begin{array}
[c]{c}%
M_{L}:\left\{
\begin{array}
[c]{c}%
x^{\prime}=\tau_{L}x+y+(\tau_{L}\xi+\eta-\xi)\ \\
y^{\prime}=-\delta_{L}x-(\delta_{L}\xi+\eta)\ \ \ \ \ \ \ \ \ \
\end{array}
\right.  \qquad \text{for } x<h-\xi\\
M_{R}:\left\{
\begin{array}
[c]{c}%
x^{\prime}=\tau_{R}x+y+(\tau_{R}\xi+\eta-\xi)\ \\
y^{\prime}=-\delta_{R}x\ -(\delta_{R}\xi+\eta)\ \ \ \ \ \ \ \ \
\end{array}
\right.  \qquad \text{for } x>h-\xi
\end{array}
\right.  \label{M}%
\end{equation}
where the prime symbol denotes the unit time advancement operator, and
$h\neq0.$

The system (\ref{M}) is topologically conjugate to\ a 2D discontinuous PWL
homogeneous map in which the functions in both partitions have the same fixed
point at the origin, $(x,y)=(0,0),$ denoted by $O.$ In fact, through the
change of coordinates $u=x+\xi,$ $v=y+\eta,$ we obtain the following map for
$X=(u,v)^{T}$, which we can rename as\ $X=(x,y)^{T}:$
\begin{equation}
T_{1}=\left\{
\begin{array}
[c]{c}%
T_{L}:X^{\prime}=J_{L}X\ \ \text{for } x<h,\ \ \ J_{L}=\left[
\begin{array}
[c]{cc}%
\tau_{L} & 1\\
-\delta_{L} & 0
\end{array}
\right]  \\ \ \\ 
T_{R}:X^{\prime}=J_{R}X\ \ \text{for } x>h,\ \ \ J_{R}=\left[
\begin{array}
[c]{cc}%
\tau_{R} & 1\\
-\delta_{R} & 0
\end{array}
\right]
\end{array}
\right.  \label{MapT}%
\end{equation}
Here, the discontinuity set is the vertical line $x=h\neq0$. Maps $M$ and
$T_{1}$ have the same dynamics, as they are topologically conjugate. 

A peculiarity of these maps, as highlighted in \cite{GRSSW-25a}, is the
emergence of a new type of attractor, referred to as a WQA. In a 2D
discontinuous map, a WQA is an attractor\footnote{An attractor $A$ of a map
$T$ is a closed invariant set with a dense orbit, for which a
neighborhood $U$ of $A$ exists such that $A=\cap_{n\geq0}T^{n}(U(A))$. In our
work, "invariant set $A$" means that it is mapped exactly into itself,
$T(A)=A.$ In other definitions, it may also be mapped strictly into itself,
$T(A)\subseteq A.$} that does not include any periodic point, thus, it is
neither an attracting cycle nor a chaotic attractor.\footnote{We refer to the
most widely used definition of chaos, that is: a map $T$ (in any dimension) is
chaotic in a closed invariant set $A$ if periodic points are dense in $A$ and
there exists an aperiodic trajectory dense in $A$ (so that there is
transitivity).} A WQA appears as the closure of quasiperiodic trajectories,
where the term "weird" refers to the rather complex and often intricate
geometric structure of these attractors. Note that if a 2D map has other
invariant sets, where it is reducible to a 1D map (e.g., a closed invariant
curve, or a set consisting of a finite number of segments), then the related
attractors are not classified as weird (although these sets may coexist with a WQA).

In \cite{GRSSW-25a}, we identified several parameter regions where map $T_{1}$
in (\ref{MapT}) with $h=-1$ has a WQA. This occurs when the fixed point $O$ is
attracting for the two linear functions, or when it is attracting only for one
function, and also when the fixed point $O$\ is repelling for both linear
functions. An example of the last case is given in Fig.1(c), showing the
coexistence of two WQAs. 

\begin{figure}
	\begin{center}
		\includegraphics[width=1\textwidth]
		{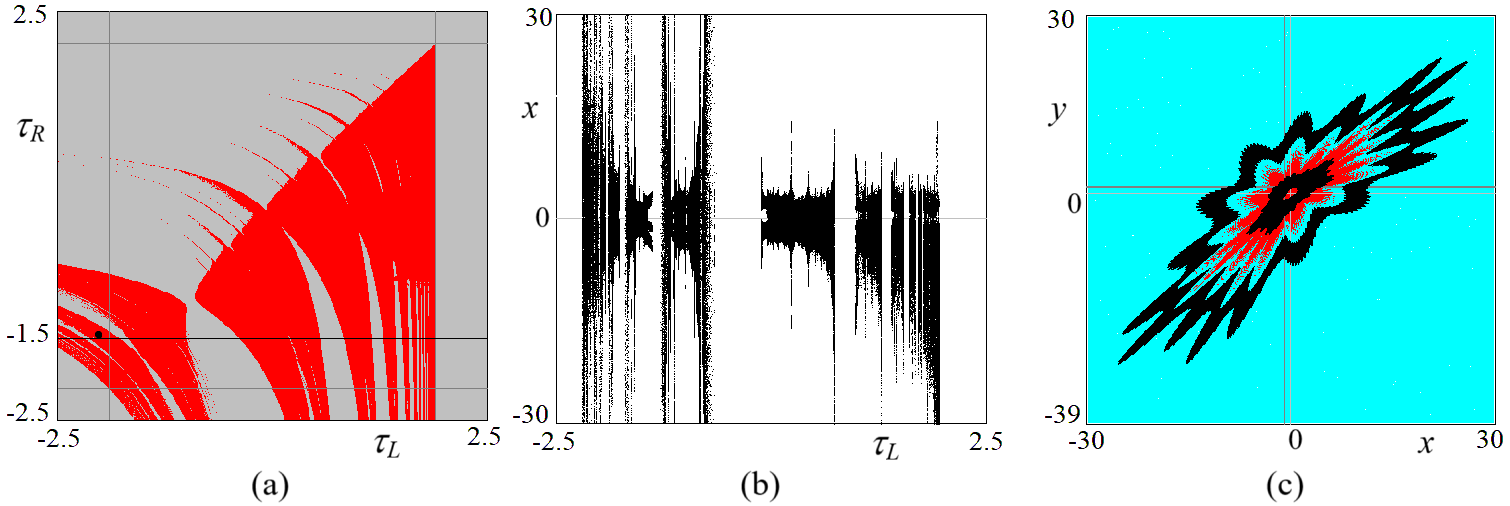} 
		\caption{\label{f1} \small{In (a), 2D bifurcation diagram in the parameter plane 
			($\tau_{L},\tau_{R})$ for map $T_{1}$ in (\ref{MapT}), with $\delta_{R}=1.11$ and
			$\delta_{L}=0.9.$ The origin is unstable for $T_{R}$, while for $T_{L}$ the
			origin is a virtual attractor in the strip between the two vertical lines, at
			$\tau_{L}=\pm(1+\delta_{L})$. In (b), 1D bifurcation diagram as a function of
			$\tau_{L}$ at $\tau_{R}=-1.5.$ In (c), phase plane at $\tau_{L}=-2$ and
			$\tau_{R}=-1.449.$ Map $T_{1}$\ has a region Z2 between $y=0.9$ 
			($=\delta_{L})$ and $y=1.11$ ($=\delta_{R})$.}}
	\end{center}
\end{figure}

In Fig.1(a), the wide red region in the $(\tau_{L},\tau_{R})$ parameter plane
of map $T_{1}$ in (\ref{MapT}) is associated with WQAs, while the gray region
is related to divergent trajectories. At the considered value of 
$\delta_{R},$ the fixed point $O$ is repelling, and the red region spreads also to the
region $\tau_{L}<-(1+\delta_{L})$, where the origin is unstable (virtual)
for function $T_{L}$ as well. In Fig.1(b), a 1D bifurcation diagram (varying
$\tau_{L}$ along the black line in (a), at $\tau_{R}=-1.5$) shows how these
attractors can mimic chaotic behaviors, although it is clear that chaos cannot
exist. In Fig.1(c), two coexisting WQAs and their basins of attraction are
shown, at the parameter point marked by a black dot in (a). At these values of
the parameters, the fixed point $O$ is repelling for the functions in both
partitions, a focus for $T_{R}$ and a saddle for $T_{L}$.

The discontinuity set in map $T_{1}$ in (\ref{MapT}) is a vertical line, with
the fixed point $O$ internal to one partition. However, as we show, the
existence of WQAs does not depend on the type of the discontinuity set. It can
be any straight line $y=mx+q$ with $q\neq0,$ or multiple straight lines, or
any curve(s) (as e.g., a circle which we use in some examples). In other
words, we consider a class of maps with any kind of discontinuity set
separating the phase space into an arbitrary finite number of partitions,
related to linear functions.

The fixed point may either be
internal to one partition or located at a border of two or more partitions. In
the first case, the map is continuous and differentiable at $O$ (so that it is
a virtual fixed point for the other functions in the definition of the map).
In the second case, the map is continuous but not differentiable at $O,$ and
at least one additional discontinuity set must exist to have a discontinuous map.

The characteristic property of the class of maps considered in this work is
the same in any dimension. It is defined as follows: 

\begin{Definition}[class of maps] We consider, for 
	$n\geq 1,$ the class of nD discontinuous PWL maps defined in a 
	finite number of partitions by linear functions with the same real fixed point.
\end{Definition}

The definition specifies a unique fixed point, which typically occurs as
hyperbolic in a generic map. Without loss of generality, we can translate it
to the origin $O$, so that the discontinuous PWL map becomes homogeneous. From
now on, we assume that the real fixed point is $O$. As mentioned earlier, it
can be internal to one partition or located on a border. In both cases, it may
be hyperbolic or nonhyperbolic for the functions defined in the respective
partitions. The properties of a linear map when the fixed point $O$ is
nonhyperbolic (e.g., possessing an eigenvalue $\lambda=1$ or $\lambda=-1,$ or
complex eigenvalues with modulus $1$) are well known. Clearly, this does not
present a problem, and the local dynamic behaviors can be easily described.

An economic application of a 2D map in which the fixed point $O$\ is always
nonhyperbolic (resulting in a segment of fixed points) can be found in
\cite{GRSSW-24}, where the existence of WQAs is also shown. Beyond economic
models, an engineering application of the class of 2D maps satisfying
Definition 1 is considered in \cite{Kollar}. 

It is important to note that WQAs may exist in other classes of piecewise
smooth maps, providing further motivation for examining their properties.
Piecewise smooth maps have gained significant attention due to their
application in various fields. In addition to economics, where they are widely
used (see \cite{Lu,Brianzoni,Day}, to cite a few), we can
mention applications in engineering (see, among others, 
\cite{DiBernardo,Simpson10,Dutta-10,Rak}), and in physics (\cite{Zhu}).

It is also worth mentioning that WQAs described in this class of maps differ
from nonchaotic attractors known as strange nonchaotic attractors, denoted as
SNAs (see, e.g., \cite{Grebogi-84,Feudel,Duan}). The main
difference is that, for maps satisfying Definition 1, there can be no
hyperbolic cycle other than the fixed point, nor any hyperbolic
tori. 

As we show in this work, WQAs have properties similar to the quasiperiodic
trajectories occurring in 1D PWL circle maps. For this reason, they may be
considered a generalization of such dynamics to the 2D phase plane. Increasing
the dimensionality naturally leads to a broader range of shapes and structures
of WQAs, which may take various and unusual forms.

In the next section, we first recall the properties of the maps satisfying
Definition 1 for $n=1$. We then examine 2D maps in our definition, proving
that bounded attractors, not associated with the fixed point $O$, are either
reducible to those of a PWL circle map (related to 1D maps) or give rise to 2D
WQA. In Section 3, we present several examples of 2D PWL maps within our
framework, defined as in (\ref{MapT}) but with different Jacobian matrices.
Section 4 explores additional examples with different types of discontinuity
sets. In Section 5, we consider the class of maps with WQAs in the $n$D space,
providing 3D examples. In Section 6, we present some open problems that need
further investigation. We summarize our findings and conclusions in Section 7.

\section{Discontinuous PWL homogeneous maps}

In Section 2.1, we first recall the properties of maps in Definition 1 for
$n=1,$ to be used in Section 2.2, where we turn to the properties of maps in
Definition 1 for $n=2.$

\subsection{1D discontinuous PWL homogeneous maps}

The case in Definition 1 for $n=1$ corresponds to the class of maps considered
in \cite{GRSSW-25b}. It is worth recalling those results, since one of our
goals is to show that a 2D map within our definition can lead to a 1D first
return map in some segment that is a function corresponding to Definition 1
with $n=1$.

In the case $n=1,$ and with only one discontinuity point, the 1D map takes the
form:
\begin{equation}
F=\left\{
\begin{array}
[c]{c}%
F_{L}:x^{\prime}=s_{L}x \quad \text{for } x<h\\
F_{R}:x^{\prime}=s_{R}x \quad \text{for } x>h
\end{array}
\ \ h\neq0\right.  
\label{MapF}%
\end{equation}
where $h\neq0$ is a scaling factor. Since $h<0$ is topologically conjugate to
$h>0,$ we set $h>0$ without loss of generality. When the slopes are positive,
bounded asymptotic dynamics distinct from the fixed point $O,$ 
can occur only if $F_{L}(h)>h>F_{R}(h),$\ leading to $s_{L}>1>s_{R}.$ In this
case, it is $F_{R}\circ F_{L}(h)=s_{R}s_{L}h=F_{L}\circ F_{R}(h),$ so that the
map is a circle homeomorphism in the invariant absorbing interval
$I=[F_{R}(h),F_{L}(h)]=[s_{R}h,s_{L}h].$ Its dynamics depend on the rotation
number $\rho$, which is well defined. It is the same for any point of interval
$I$ (see \cite{deMelo}). If $\rho$ is rational, then interval $I$ is filled
with periodic points (of the same period). If $\rho$ is irrational, then $I$
is filled with quasiperiodic orbits dense in the interval. The generic case is
an irrational rotation, since a rational rotation is associated with a set of
zero Lebesgue measure in the $(s_{R},s_{L})$ parameter plane of the possible
slopes. This occurs if and only if integers $p$ and $q$ exist such that
$s_{L}^{p}s_{R}^{q}=1.$ Let us assume that $p$ and $q$ are the smallest
integers. Then the nonhyperbolic cycles have period $n=p+q$ and rotation
number $\rho=p/n$ (or $(1-\rho)=q/n$).\footnote{It is worth noting that the
map known as "rigid rotation", that is 
$R:$ $x\rightarrow x+\rho (\operatorname{mod}1),$ is not included in our definition 
(since the fixed point is at infinity). However, this circle homeomorphism with rotation number
$\rho$ is topologically conjugate to map $F$ with positive slopes and the same
rotation number.} We also mention that in the class of 1D PWL\ Lorenz maps
(with one discontinuity point), the case of a circle map denotes the
transition from regular dynamics (in a gap map, where chaos cannot occur) to
chaotic dynamics (in an overlapping map), see \cite{Rand,Berry,Avrutin}.

The dynamics of the generic map for $n=1$ in\ Definition 1 is described in the
following

\begin{Theorem}[from \cite{GRSSW-25b}] Let $G$ 
be a 1D discontinuous PWL homogeneous map as in Definition 1. Then: 
\begin{itemize}[]

\item[(1)] A hyperbolic cycle different from the fixed point $O$ 
cannot exist (and thus, a chaotic set cannot exist).

\item[(2)] The only possible bounded invariant sets of map $G,$ 
different from those related to the fixed point $O$ 
(whether hyperbolic or nonhyperbolic), are those occurring in a PWL circle
map. These consist of intervals densely filled with nonhyperbolic cycles or
quasiperiodic orbits. Coexistence is possible.

\item[(3)] Quasiperiodic orbits lead to (weak) sensitivity to initial
conditions.

\item[(4)] The Lyapunov exponent is zero.
\end{itemize}
\end{Theorem}

\subsection{2D discontinuous PWL homogeneous maps}

In this section, we characterize the attracting sets that may occur in 2D PWL
maps as defined in Definition 1. However, the properties described in the
following lemma hold in any dimension: 

\begin{Lemma} Let $T$ be an $nD$ map as in
Definition 1. Then: 
\begin{itemize}[]

\item[(i)] A hyperbolic cycle different from the fixed point $O$ 
cannot exist (and thus, no chaotic set in $nD$ can exist).

\item[(ii)] Segments of straight lines through the fixed point $O$ 
are mapped into segments of straight lines through $O$. 

\item[(iii)] Any composition of the linear functions defining map 
$T$ preserves properties (i) and (ii).
\end{itemize}
\end{Lemma}

Proof. 
(i) Let us assume that a hyperbolic $k-$cycle\ different from the fixed point
$O$ exists, with periodic points in the partitions labelled as usual via the
ordered symbolic sequence, say $\sigma$ (with $k$ symbols). Then such a cycle
must lead to $k$ hyperbolic fixed points of map $T^{k}$, distinct from $O$.
But this is not possible, because the composition of linear homogeneous
functions is always a linear homogeneous function with $O$ as the unique fixed
point, provided that it is hyperbolic, as assumed. Therefore, no hyperbolic
cycle other than $O$ can exist, and a chaotic set in $nD$ is impossible.

(ii) This is an immediate consequence of the considered class of PWL
homogeneous maps, since any linear function maps segments of straight lines
into segments of straight lines. Thus, in particular, this holds for straight
lines through $O$, fixed point for all functions defining map $T$.

(iii) Also this property is immediate, since any composition of linear
homogeneous functions is a function in the same class.$\square$ 
\medskip 

As mentioned in the Introduction, it is possible for the dynamics of map $T$
to involve a segment $\tau$ of some particular straight line. Let $r$ be a
straight line through the fixed point $O,$ $y=mx$. Consider a point
$(x,mx)\in\tau\subset r$ within a partition of $T.$ Suppose that after a
finite number of iterations, say $k$, another point 
$(x^{\prime},mx^{\prime})$ on $\tau\subset r$ is obtained. This implies:%
\begin{equation}
x^{\prime}\left[
\begin{array}
[c]{c}%
1\\
m
\end{array}
\right]  =xA_{k}\left[
\begin{array}
[c]{c}%
1\\
m
\end{array}
\right]  ,\ \ A_{k}=[J_{j_{k}}...J_{j_{1}}] 
\label{p1}%
\end{equation}
where $A_{k}$ is the product of the $k$ Jacobian matrices of the functions
applied during the trajectory. This equation shows that $r$ is an eigenvector
of matrix $A_{k}$, associated with the eigenvalue 
$\lambda=\frac{x^{\prime}}{x}.$

In particular, when matrix $A_{k}$ has the characteristic
polynomial  $\mathcal{P}(\lambda)=\lambda^{2}-Tr(A_{k})\lambda+\det(A_{k})$
satisfying $\mathcal{P}(1)=0,$ then this implies that $r$ is an eigenvector of
matrix $A_{k}$ associated with the eigenvalue $\lambda=1$ (since 
$x^{\prime}=x$). This may lead to nongeneric cases, characterized by segments filled
with nonhyperbolic cycles (nonhyperbolic fixed points of $T^{k}$). The
symbolic sequence of the cycles corresponds to the ordered Jacobian matrices
involved along the trajectory, say $\sigma=j_{1}....j_{k}$, where $j_{i}$
identifies a partition, assuming that the admissibility conditions related to
the partitions are satisfied. If we consider the first return map on the
appropriate segment of the eigenvector, we obtain the identity function. This
holds cyclically, on $k$ segments, none of which can intersect a discontinuity
set (as that would lead to different symbolic sequences for periodic points of
the same segment). However, a discontinuity point and other critical points
(images of the discontinuity point) form the boundaries of the invariant
segments. 

In our previous work \cite{GRSSW-25a}, we identified particular regions in the
parameter space of map $T_{1}$ in (\ref{MapT}), where WQAs may exist and
persist under parameter perturbation. We have shown that a mechanism that
leads to their appearance may be linked to particular cyclical invariant
segments, or half-lines, filled with nonhyperbolic cycles. More precisely,
particular attracting sets may be associated with these sets in the parameter
space, occurring when the characteristic polynomial of some matrix $A_{k}$
satisfies $\mathcal{P}(1)=0$ and each of the cyclical segments belongs to the
proper partition (admissibility condition). Such an attracting set,
structurally unstable, may play a key role in the appearance of WQAs, when
parameters are perturbed.

A mechanism leading to the appearance of WQAs similar to the one described
in \cite{GRSSW-25a} will be shown in the following sections, including
examples in maps with different kinds of discontinuity sets. 

Another particular case where a straight line plays a role in the dynamics of
map $T$ arises, for example, when the Jacobian matrix in a given partition has
an eigenvalue $\lambda=0$ (specific examples will be provided in later
sections). In such cases, it may be possible to define the first return map on
a segment of the eigenvector associated with the eigenvalue 
$\lambda=0.$

In a generic case, to construct the first return map on a segment $\tau$ of a
straight line $r$ through the fixed point $O$ (say, $y=mx$)$,$ each
point $(x,mx)$ is iterated and when the trajectory returns, for the first
time, to the same segment, say $(x^{\prime},mx^{\prime})$ after $k_{1}$
iterations, it satisfies: 
\begin{equation}
x^{\prime}\left[
\begin{array}
[c]{c}%
1\\
m
\end{array}
\right]  =xA_{k_{1}}\left[
\begin{array}
[c]{c}%
1\\
m
\end{array}
\right]  ,\ \ A_{k_{1}}=[J_{j_{k_{1}}}...J_{j_{1}}] 
\label{first}%
\end{equation}
where, as above, $A_{k_{1}}$ is the product of the $k_{1}$ Jacobian matrices
of the functions applied during the trajectory. Then we assign to $x$ the
value $x^{\prime}$, which is obtained via eigenvalue 
$\lambda_{k_{1}}=\frac{x^{\prime}}{x},$ defining $F_{r,k_{1}}(x)=\lambda_{k_{1}}x.$ That is,
$x^{\prime}=\lambda_{k_{1}}x.$ Considering the symbolic sequences of the
trajectory in the usual way, using the partitions involved, the first return
follows a fixed symbolic sequence, say $\sigma=j_{1}....j_{k_{1}}.$ Due to the
linearity, the definition remains valid for an interval of points along the
segment (i.e., all points have the same symbolic sequence, and thus the same
return map), up to a point, say $(c,mc)$, whose trajectory merges with a
discontinuity point of map $T$.

If the first return exists on the same segment also for $x>c,$ it means that a
second integer $k_{2}$ exists such that after $k_{2}$ iterations a point
$(x^{\prime},mx^{\prime})$ belongs to the segment:
\[
x^{\prime}\left[
\begin{array}
[c]{c}%
1\\
m
\end{array}
\right]  =xA_{k_{2}}\left[
\begin{array}
[c]{c}%
1\\
m
\end{array}
\right]  ,\ \ A_{k_{2}}=[J_{j_{k_{2}}}...J_{j_{1}}]
\]
and the first return leads to a different function, say 
$F_{r,k_{2}}(x)=\lambda_{k_{2}}x$, related to a different symbolic sequence. In such a
case, the first return map has two partitions, each defined by a linear
homogenous function (with the same fixed point in the origin), given by
$x^{\prime}=\lambda_{k_{1}}x$ for $x<c$ and $x^{\prime}=\lambda_{k_{2}}x$ for
$x>c.$ That is, we have a map as map $F$ in (\ref{MapF}).

Let us now prove the following 

\begin{Theorem} Let $T$ be a 2D map as in Definition 1.
Let $r$ be a straight line through the fixed point $O$ 
($y=mx$), such that the first return of map $T$ on a segment
of $r$ leads to a 1D map with a finite number of discontinuity
points. Then the related dynamics of map $T$ in the phase plane are
either nonhyperbolic cycles (with eigenvalue 1) or quasiperiodic trajectories
densely filling some segments.
\end{Theorem}

Proof. 
The proof consists in showing that the first return is a 1D PWL discontinuous
map defined by homogeneous functions (with fixed point $O$). Consider a point
$(x,mx)\in r,$ applying the functions up to its return on $r$, a point
$x^{\prime}$ as in (\ref{first}) is obtained. By continuity, the map remains
$x^{\prime}=\lambda_{k_{1}}x,$ where $\lambda_{k_{1}}$ is an eigenvalue of
$A_{k_{1}},$ up to a point $\xi_{1}$ whose trajectory merges with a
discontinuity point of map $T$. For $x>\xi_{1},$ the first return changes
definition, symbolic sequence, and eigenvalue, leading to a different linear
function $x^{\prime}=\lambda_{k_{2}}x,$ where $\lambda_{k_{2}}$ is an
eigenvalue of $A_{k_{2}}.$ This process continues. If a finite number of
partitions (or discontinuity points) and matrices $A_{j},$ $j=k_{1},...,k_{m}$
form the first return map on $r$, then the map is a 1D map $G(x)$ that
satisfies Definition 1 for $n=1.$ It follows that the dynamics of map $G$ (and
map $T$) are necessarily as those described in Theorem 1, associated with a
PWL circle map.$\square$ \medskip 

\begin{figure}[h]
	\begin{center}
		\includegraphics[width=0.6\textwidth]
		{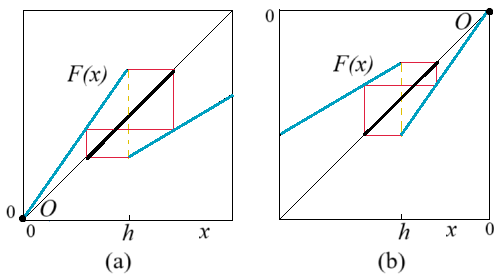} 
		\caption{\label{f2} \small{Qualitative representation of the first return map on a segment. In (a),
			at the right side of the fixed point $O$. In (b), at the left side of the
			fixed point $O$.}}
	\end{center}
\end{figure}%

Thus, when an invariant segment exists on a straight line through $O$, 
the first return map can be reduced to
a 1D map satisfying Definition 1 for $n=1$. That is, a map $x^{\prime}=G(x)$
with many discontinuity points, and linear branches with the fixed point $O$.
However, a suitable segment exists where the first return is the simplest one,
with only one discontinuity point, corresponding to the form of map $F$ in
(\ref{MapF}) (see \cite{GRSSW-25b}). 

Moreover, in the case of only one discontinuity point, when the first return
is in a segment on the right side of the fixed point $O$, starting from a
point $x>0$ left boundary of the interested interval$,$ it must be 
$x^{\prime}>x>0$ so that it is $\lambda=\frac{x^{\prime}}{x}>1.$ If there are bounded
dynamics, this also implies, as commented for map $F$ in (\ref{MapF}) for
$h>0$, that the second branch must necessarily be with $\lambda<1$, as in
Fig.2(a). If the segment is at the left side of fixed point $O$, starting from
the left boundary of the interested interval, it must be $x<x^{\prime}<0$ so
that it is $0<\lambda=\frac{x^{\prime}}{x}<1.$ This also implies that if there
are bounded dynamics (different from the fixed point $O$), the second branch
must be with $\lambda>1$, as in Fig.2(b) (case $h<0$, conjugate to the one
with $h>0$, as commented for map $F$ in (\ref{MapF})). Several examples
illustrating these cases are presented in the following sections. 

We give now the main result for 2D maps:

\begin{Theorem} Let $T$ be a 2D map as in Definition 1. Then: 
\begin{itemize}[] 
\item[(j)] A bounded $\omega$-limit set $\mathcal{A}$ 
different from the fixed point $O$ and from local invariant sets
associated with $O$ when it is nonhyperbolic, 
can only be one of the following kind:  
\begin{itemize}[] 
\item[(ja)] it is a nonhyperbolic $k$-cycle, $k\geq2$, belonging to $k$ segments 
not intersecting any border, filled with $k$-periodic points (all cycles have the same symbolic sequence); 

\item[(jb)] it belongs to an invariant set on which the dynamics are reducible to a discontinuous 1D map; 

\item[(jc)] it is a weird quasiperiodic attractor (WQA).
\end{itemize} 
\item[(jj)] When $\mathcal{A}$ does not consist in segments filled with cycles, 
then map $T$ exhibits (weak) sensitivity to initial conditions in $\mathcal{A}$. 
\end{itemize}
\end{Theorem}

Proof.  
From Lemma 1, it follows that the map cannot have hyperbolic cycles
different from the fixed point $O$. If a $k$-cycle exists, $k\geq2$, 
it can only be nonhyperbolic.

Nonhyperbolic $k$-cycles may be dense in $k$ intervals, belonging to
eigenvectors of components of the map $T^{k}$, leading to (ja).

When a suitable restriction of $T$, or a first return map, can be defined in
some interval, then the dynamics of $T$ are reducible to those of a discontinuous 1D map,
leading to case (jb). Moreover, when a first return map exists in some
interval on a straight line through $O$, different from the identity function, 
it follows from Theorem 2 that it is a 1D map $G(x)$ as in Definition 1 for 
$n=1$, and by Theorem 1, its dynamics are reduced to those of a PWL circle map. 
In this case, $\mathcal{A}$ consists in a finite number of segments
filled with quasiperiodic orbits or filled with periodic points.

(jc) If the dynamics are not reducible to those of a 1D map then, recalling
that $\mathcal{A}$ is a closed invariant set in 2D, given by the limit points of a
trajectory, we can only have quasiperiodic trajectories dense in 
$\mathcal{A}$. This follows because the trajectory of a point 
$(x_{0},y_{0})\in \mathcal{A}$ can never come back to the same point (that would lead to a
cycle). Thus, it can neither be a chaotic set (since it does not include any
cycle). Then $\mathcal{A}$ is a weird quasiperiodic attractor in 2D. 

(jj) Sensitivity to initial conditions follows immediately because, for any
two nearby points within the invariant set $\mathcal{A}$, in a finite number
of iterations their trajectories will be on opposite sides of a discontinuity.
Once this occurs, the points are mapped by different functions and they leave each
other's neighborhood. The sensitivity is weak in case (jb), since it is
associated with a 1D map in a 2D phase plane. Due to the absence of chaotic
behavior, sensitivity is weak also in case (jc).$\square$ 
\medskip 

Clearly, in the same map $T$ different kinds of invariant sets can coexist. 
Examples of coexisting attracting sets and WQAs in the class of 2D maps in
Definition 1 have been shown in \cite{GRSSW-24,GRSSW-25a,GRSSW-25c}. 
In the following sections, we provide several examples in
PWL maps having different types of Jacobian matrices and various discontinuity
sets. 

It is important to remark that in the parameter space of the considered 2D
maps, the generic attractor different from the fixed point, when existing, is
a WQA. In fact, the attracting sets related to (ja) and (jb) in Theorem 3 are
not structurally stable. For example, they occur when a suitable composition of the
functions exists, leading to particular segments filled with nonhyperbolic
cycles (associated with an eigenvalue equal to $1$), satisfying admissibility
conditions. Or when some peculiar properties (as the case of $0$ as an
eigenvalue) may lead to an attracting set belonging to a finite number of
segments in the phase plane. As we shall see in the examples, these particular
cases may coexist with a WQA. We shall also show that even if half a plane is
mapped into a straight line, the dynamics may lead to a WQA. These cases
associated with (ja) and (jb) in Theorem 3 may exist, but are not persistent.
They occur for a set of zero Lebesgue measure in the parameter space, and a
small perturbation of the parameters in general leads to the appearance of a WQA.

As for the 2D bifurcation diagram in Fig.1(a), and for all the 2D bifurcation
diagrams presented in the next sections (numerical computations obtained
starting with a given initial condition in the phase plane), gray color
indicates divergence of the trajectory, green color indicates convergence to
the fixed point $O$, while red color corresponds to convergence to an
attractor, which may be either a WQA (which is a generic case (jc) in Theorem
3) or another attracting set (nongeneric cases (ja) and (jb) in Theorem 3). It
is worth mentioning that also the 2D bifurcation diagrams reported in the
examples in the next sections with a fixed value $\delta_{L}=0$ (showing large
red regions) may be related to attracting sets that are not structurally
stable, or to WQAs.

\section{2D PWL maps with discontinuity set $x=-1$}

\subsection{PWL map $T_{1}$}

Let us consider the map in (\ref{MapT}) with $h=-1$. As in \cite{GRSSW-25a},
the discontinuity line is referred to as a critical line, as well as its
images by the two linear functions, given by
\begin{equation}
T_{L/R}(x=-1):\ \ y=\delta_{L/R}%
\end{equation}
Assuming that no eigenvalue is equal to zero (i.e. $\delta_{L/R}\neq0$), each
half-plane bounded by $x=-1$ is mapped by the linear function $T_{L/R}$ into a
half-plane bounded by $y=\delta_{L/R}$. The relative positions of the
half-planes determine the classification of the kind of map. In particular,
for $\delta_{L}\neq\delta_{R}$, the map may be uniquely invertible or
noninvertible. The strip bounded by $y=\delta_{L/R}$ may be a so-called
region $Z_{0},$ whose points have no rank-1 preimage, or a region $Z_{2},$
whose points have two different rank-1 preimages.

In the example shown in Fig.1(c), the 2D map $T_{1}$ has a region $Z_{2}$, and
points of the WQA also belong to that strip (however, it is possible that for
such points of the attractor, only one rank-1 preimage belongs to the
attractor). One more example with a region $Z_{2}$ is shown in Fig.3(c).
Differently, in the example in Fig.3(b), the map has a region $Z_{0}$. It is
important to note that when a region $Z_{0}$ exists, no point of an invariant
set can belong to that region, nor to any of its image of any rank. Thus, in
such a case, we can conclude that any invariant set $\mathcal{A}$ of the map
(in particular a WQA) must belong to the region 
$\mathbb{R}^{2}\backslash \underset{n\geq0}{\cup}T^{n}(Z_{0})$. 

\begin{figure}
	\begin{center}
		\includegraphics[width=1\textwidth]
		{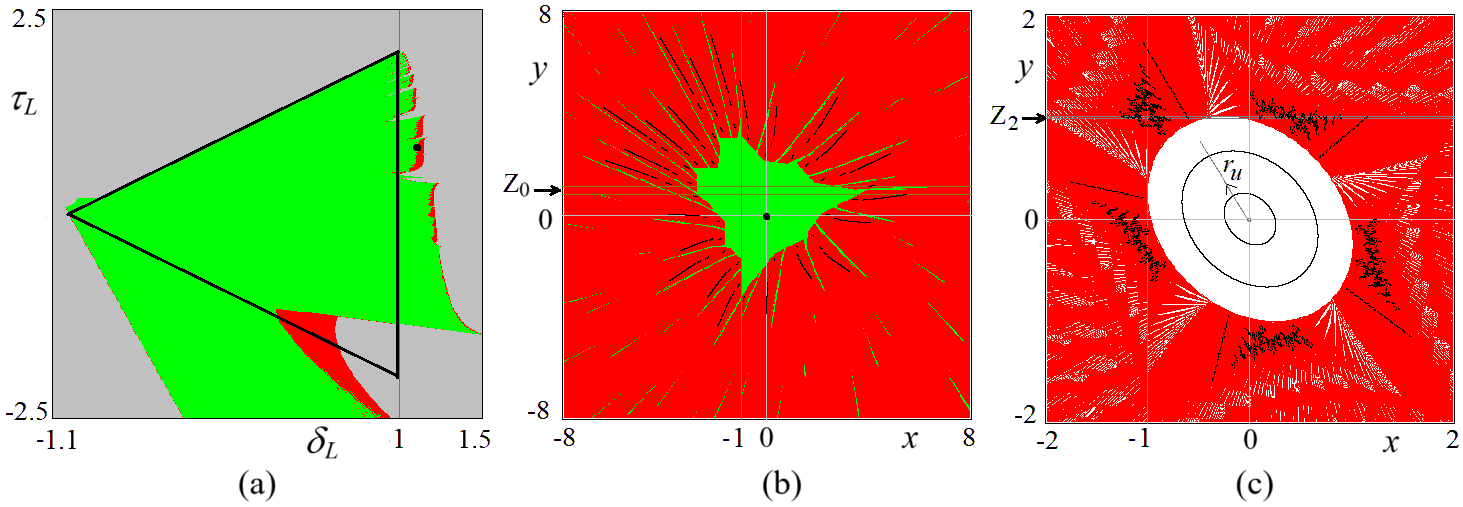} 
		\caption{\label{f3} \small{In (a), 2D bifurcation diagram in the ($\delta_{L},\tau_{L})$ parameter
			plane\ for map $T_{1}$ in (\ref{MapT}), with $\delta_{R}=0.8$ and 
			$\tau_{R}=-0.5.$ The black lines denote the standard stability triangle for the
			linear map $T_{L},$ bounded by segments of the lines of equation 
			$\tau_{L}=\pm(1+\delta_{L})$ and $\delta_{L}=1.$ In (b), phase plane at 
			$\delta_{L}=1.1$ and $\tau_{L}=1$ (black dot in (a)). Map $T_{1}$ has a region
			$Z_{0}$ between $y=\delta_{R}$ and $y=\delta_{L}.$ In (c), phase plane at
			$\delta_{R}=1,$ $\tau_{R}=-0.5,$ $\delta_{L}=0.98$ and $\tau_{L}=0.8.$ Map
			$T_{1}$\ has a region $Z_{2}$ between $y=\delta_{L}$ and $y=\delta_{R}$. The
			fixed point $O$ is a center, with an irrational rotation number; the large
			white region in the $R$ partition is filled with invariant ellipses, bounded
			by the ellipse tangent to the discontinuity line and all its images.}}
	\end{center}
\end{figure}

In Fig.3(a), the parameters of the function in the right partition are set
such that the origin is an attracting focus for map $T_{1}$. The 2D
bifurcation diagram illustrates the stability triangle of function $T_{L}.$
The black dot is in a region where the fixed point $O$ is a repelling focus
for function $T_{L},$ and a weird quasiperiodic attractor coexists with the
attracting origin. These two attractors are shown in Fig.3(b) in black, along
with the related basins of attraction (in green, the basin of the fixed point
$O$, in red the basin of the WQA).

It is important to note that in these PWL maps, the boundaries of the basins
of attraction include segments of the discontinuity set and related preimages
of any rank \cite{Mira-96}.

In the previous section, we mentioned that the fixed point $O$ may be
nonhyperbolic. An example is shown in Fig.3(c), with $\delta_{R}=1;$ the fixed
point $O$ is a center in the right partition. In our example, the rotation
number is irrational, meaning that in the $R$ partition there exists an
invariant region bounded by an ellipse tangent to the discontinuity line and
its images. This region is filled with ellipses, on which the trajectories are
quasiperiodic (see \cite{SG-2008}). However, a WQA coexists. Furthermore,
since in this case the map has a region $Z_{2}$ that includes a portion of the
invariant region in the right partition, this region has preimages (shown in
white in Fig.3(c)) that belong to the stable set, or basin of attraction, of
the invariant region.

In the example shown in Fig.3(c), the mechanism of appearance of the WQA is
similar to that described in \cite{GRSSW-25a}. In fact, considering the
composite map $T=T_{L}\circ(T_{R})^{4}$ $($i.e., $T(X)=T_{L}(T_{R}^{4}(X))$),
the origin is a (virtual) saddle of $T,$ with eigenvalues 
$\lambda_{1}\simeq0.94$ and $\lambda_{2}\simeq1.042.$ For $0<\lambda_{2}<1,$ the invariant
polygon attracts all the points of the phase plane. For $\lambda_{2}=1,$ there
exist five invariant segments (bounded on both sides), filled with 5-cycles,
fixed points of map $T$ (as commented in (\ref{p1}), and symbolic sequence
$\sigma=RRRRL$). For $\lambda_{2}>1,$ the trajectories on the eigenvector, say
$r_{u}$, become repelling. A segment of $r_{u}$ in Fig.3(c) shows the
direction, although it applies only outside the invariant polygon. Let $P$ be
the intersection point of eigenvector $r_{u}$ with the discontinuity line, and
$P_{-1}$ its rank-1 preimage (which lies within the $R$ partition). Function
$T$ maps segment $P_{-1}P$ into a segment $PP_{1}$ (along the eigenvector).
This segment now belongs to the $L$ partition, where a different function is
applied. Then the iterates of this segment converge to a WQA. We can say that
the WQA is the $\omega$-limit set of $T_{1}^{n}(PP_{1}),$ for 
$n\rightarrow \infty.$

\subsection{PWL map $T_{1}$ with $\delta_{L}=0$}

In this class of maps, certain peculiar cases are worth describing. In
particular, when one of the two linear functions has an eigenvalue equal to
zero, a whole region of the phase plane (in this case, a half-plane) is mapped
into the corresponding image of the discontinuity line, that is an eigenvector
of the Jacobian matrix in that partition.

When this occurs, it may be that the dynamics of the 2D map can be
investigated via a 1D first return map on that line. This is because any
existing bounded invariant set, different from those associated with the fixed
point $O$, must necessarily have points on that critical line.

In Fig.4(a), we present a 2D bifurcation diagram of map $T_{1}$ with
$\delta_{L}=0.$ In Fig.4(b), we consider a specific parameter point (marked by
a black dot in Fig.4(a)), at which an attractor consisting of a finite number
of segments coexists with the attracting fixed point $O$ (whose basin is shown
in green). The first return map in the segment of the attractor that lies
along the critical line $y=0$ is shown in Fig.4(c), confirming that the
dynamics of the attractor (distinct from the fixed point $O$) are those of a
PWL circle map. 

\begin{figure}[h]
	\begin{center}
		\includegraphics[width=1\textwidth]
		{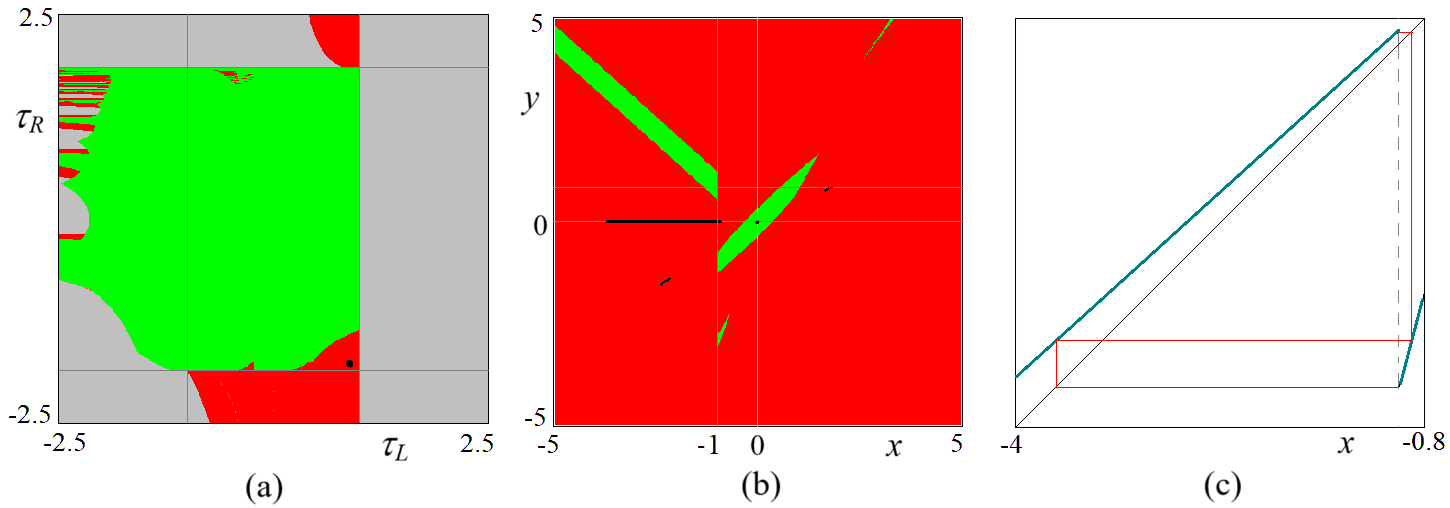} 
		\caption{\label{f4} \small{ In (a), 2D bifurcation diagram in the ($\tau_{L},\tau_{R})$ parameter
			plane for map $T_{1}$ in (\ref{MapT}), with $\delta_{R}=0.85$ and 
			$\delta_{L}=0.$ The left partition is mapped by map $T_{1}$ onto the critical line
			$y=0$. In (b), phase plane at $\tau_{L}=0.9$ and $\tau_{R}=-1.8$ (black dot in
			(a)), the attracting fixed point $O$ coexists with another attractor having a
			segment on $y=0$. In (c), first return map in the segment of the attractor
			belonging to the line $y=0,$ that is a PWL circle map.}}
	\end{center}
\end{figure}

When the origin is unstable, another attractor may exist for a wide range of
parameter values, as evidenced in the following example. In Fig.5(a), the
parameters of the function in the right partition are set so that the origin
is a repelling focus. The 2D bifurcation diagram shows the stability triangle
of function $T_{L}.$ In particular, the black dot belongs to the region in
which the origin is attracting for function $T_{L}$ with $\delta_{L}=0.$ The
corresponding attracting set is shown in Fig.5(b). Also here, the first return
map on a segment on the critical line $y=0$ confirms that the dynamics of the
attractor are those of a PWL circle map.

\begin{figure}
	\begin{center}
		\includegraphics[width=1\textwidth]
		{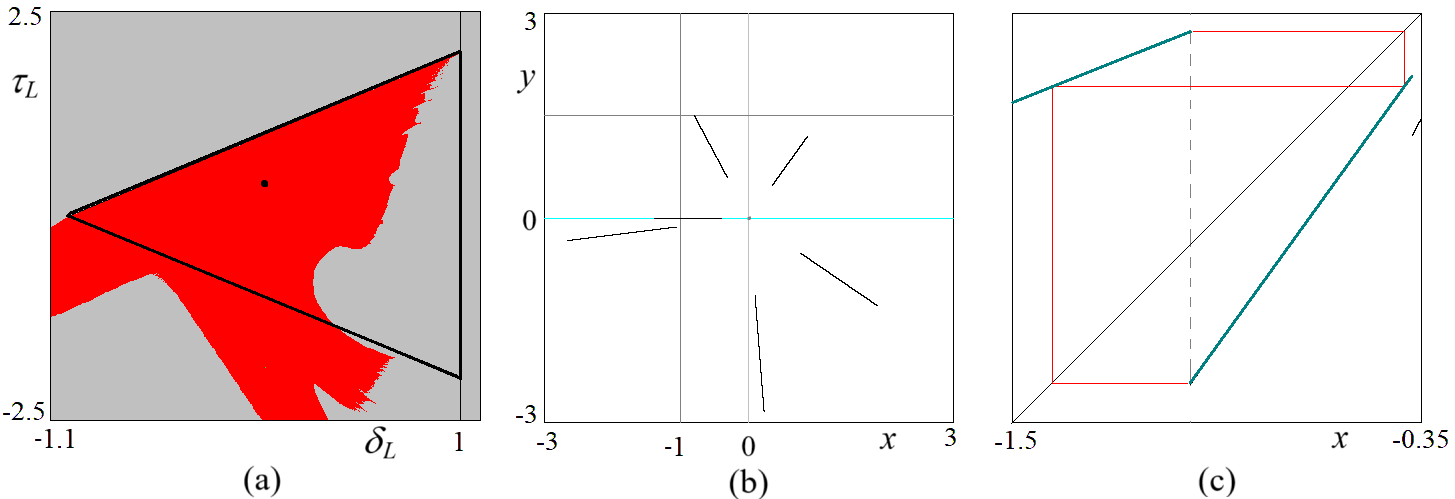} 
		\caption{\label{f5} \small{In (a), 2D bifurcation diagram in the $(\delta_{L},\tau_{L})$ parameter
			plane for map $T_{1}$ in (\ref{MapT}), with $\delta_{R}=1.5$ and 
			$\tau_{R}=0.8.$ The stability triangle of the virtual fixed point is highlighted.
			At the black dot, $(\delta_{L},\tau_{L})=(0,0.4),$ the left partition is
			mapped by map $T_{1}$\ onto the critical line $y=0$. The attractor existing in
			the phase plane is shown in (b). In (c), first return map in the segment of
			the attractor belonging to line $y=0,$ that is a PWL circle map.}}
	\end{center}
\end{figure}

Even when the attracting set consists of several intervals along the critical
line $y=0$, the first return map can be defined on a single interval,
resulting in a PWL circle map, as illustrated in the example of Fig.6(a,b). 

\begin{figure}
	\begin{center}
		\includegraphics[width=0.8\textwidth]
		{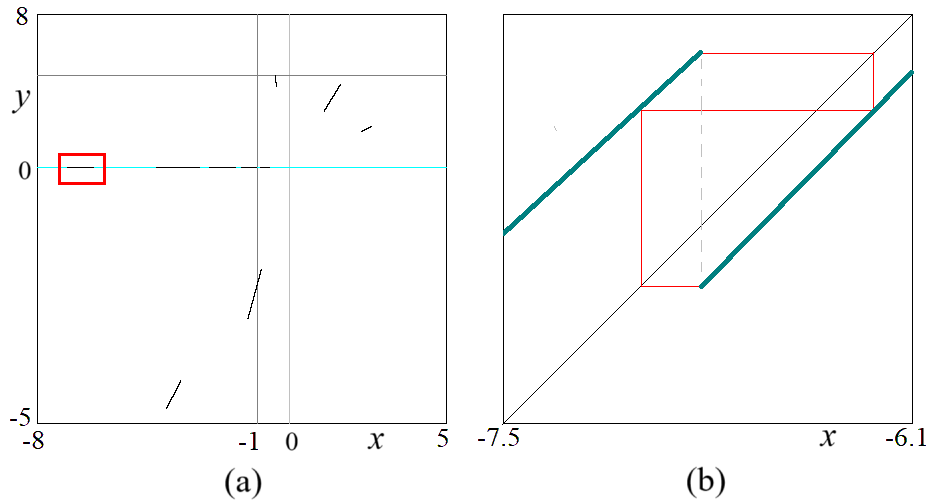} 
		\caption{\label{f6} \small{In (a), quasiperiodic attractor of map $T_{1}$ in (\ref{MapT}) in the
			phase plane at $\delta_{L}=0,$ $\tau_{L}=0.6,$ $\delta_{R}=1.8$ and 
			$\tau_{R}=-1.8.$ In (b), the first return map in the segment of the attractor
			belonging to the line $y=0,$ highlighted in (a) by a red rectangle, showing a
			PWL circle map.}}
	\end{center}
\end{figure}

However, the condition $\delta_{L}=0$ (which causes the left partition to be
mapped onto the critical line $y=0$) does not necessarily imply that an
attractor consists of a finite number of segments. The example in Fig.7
suggests that this is not always the case. In Fig.7(a), the 2D bifurcation
diagram in the parameter plane $(\tau_{L},\tau_{R})$ is shown for fixed
parameters with $\delta_{R}=1.01$ (so that $O$ is unstable) and 
$\delta_{L}=0$. The attracting set corresponding to the black dot in Fig.7(a), at
which the origin is a repelling focus, is illustrated in Fig.7(b). In this
case, we are not able to define a suitable first return map on the critical
line $y=0$, and the attracting set appears to be a WQA. 

\begin{figure}
	\begin{center}
		\includegraphics[width=1\textwidth]
		{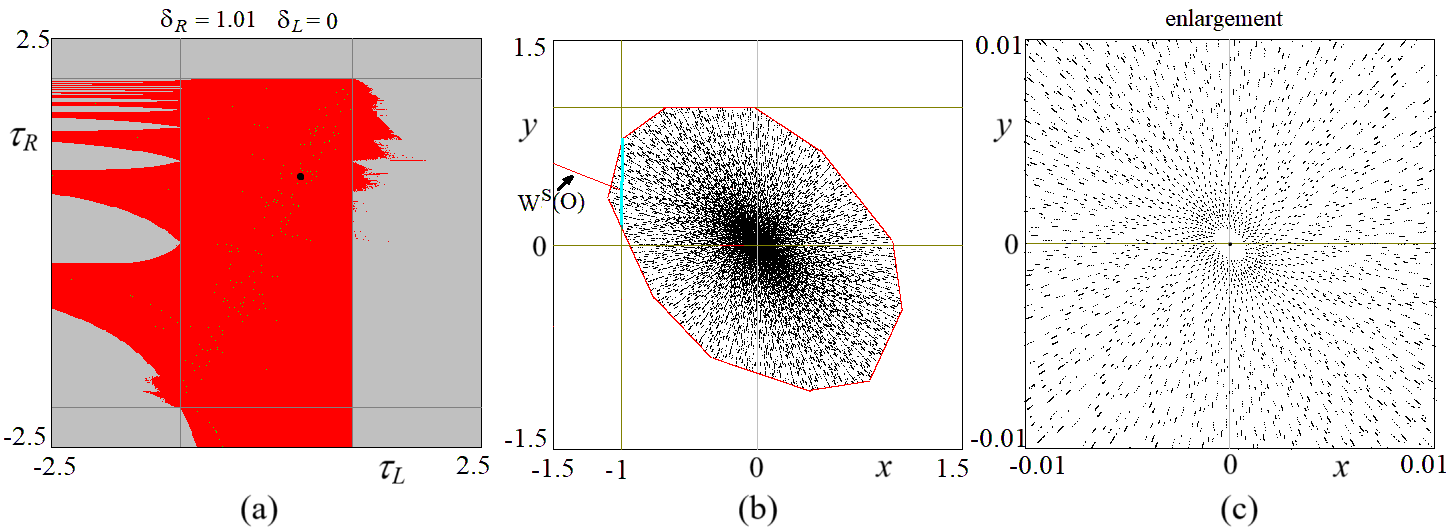} 
		\caption{\label{f7} \small{In (a), 2D bifurcation diagram in the $(\tau_{L},\tau_{R})$ parameter
			plane for map $T_{1}$ in (\ref{MapT}), with $\delta_{R}=1.01$ and 
			$\delta_{L}=0.$ The left partition is mapped by map $T_{1}$ onto the critical line
			$y=0$. The attractor existing in the phase plane at the black dot, 
			$(\tau_{L},\tau_{R})=(0.4,0.8)$, is shown in (b). In (c), the enlargement of a
			neighborhood of the fixed point $O$, repelling focus.}}
	\end{center}
\end{figure}

In this case, the origin is a hyperbolic fixed point, a repelling focus for
map $T_{R}$. An invariant polygon can be defined, including the WQA, via a
finite number of images of the segment of the discontinuity set (see
\cite{Mira-96}), colored in azure in Fig.7(b). In Fig.7(b), it is shown in red
the polygon and stable set of $O,$ $W^{S}(O).$ In one iteration, that
half-line in the $L$ partition is mapped by the map into the fixed point $O$.
It appears that no point of the WQA belongs to the hlalf-line $W^{S}(O)$. The
enlargement in a small neighborhood of the fixed point $O$ in Fig.7(c) shows
that there is a hole, a neighborhood of the origin, without points of the WQA.
This means that there are no points of the WQA close to the stable set of $O$
(half-line in the $L$ partition, shown in red in Fig.7(b)).

Clearly, the same example with $\delta_{L}=-0.01$ also leads to a WQA, shown
in Fig.8(a), and it is interesting to see how the points now fill the existing
invariant polygon (although not densely). The enlargement in Fig.8(b) shows
that here, too, there are no points of the WQA in a suitable neighbohood of
the fixed point $O$. 

\begin{figure}
	\begin{center}
		\includegraphics[width=0.8\textwidth]
		{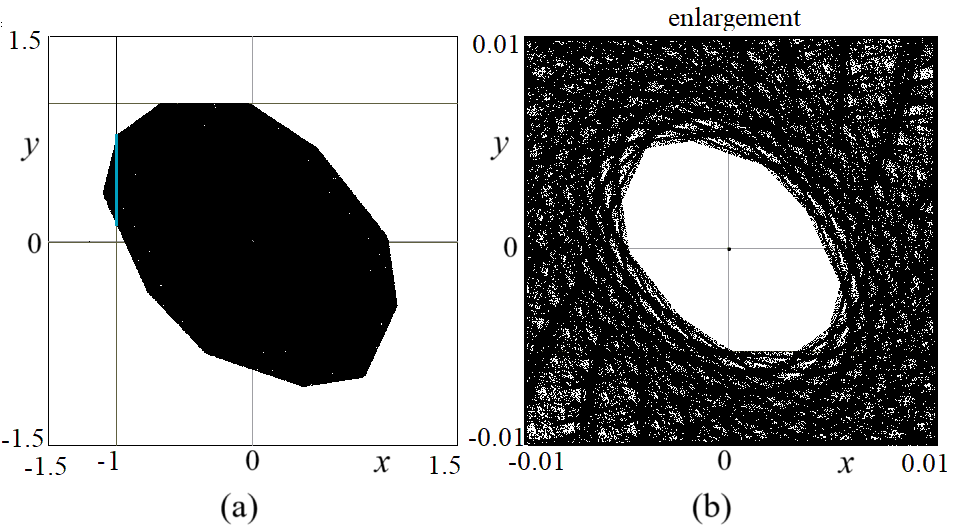} 
		\caption{\label{f8} \small{In (a), WQA in the phase plane for map $T_{1}$ in (\ref{MapT}) at
			$\tau_{R}=0.8$, $\delta_{R}=1.01$, $\tau_{L}=0.4$ and $\delta_{L}=-0.01.$ In
			(b), the enlargement of a neighborhood of the fixed point $O$, repelling
			focus.}}
	\end{center}
\end{figure}

The existence of a WQA in this case can be explained as follows. The half-line
stable set of $O$ (in the previous example, in Fig.7(b), with slope $-0.4$) is
now an eigenvector of $T_{L}$ with slope $-0.4236$ (corresponding to the
eigenvalue $\lambda=0.0236)$ whose points tend toward the virtual fixed point
$O$. In Fig.8(a), this eigenvector of $T_{L}$ intersects the discontinuity
line at a point $P$, inside the invariant region. Point $P$ also has a rank-1
preimage via the inverse $T_{L}^{-1},$ denoted by $P_{-1}$, which belongs to
the eigenvector within the left\ partition. Consequently, $T_{L}$ maps segment
$P_{-1}P$ into segment $PP_{1}$ along the eigenvector. However, segment
$PP_{1}$ now belongs to the right partition, where the right function applies,
and in a finite number of iterations, the points are mapped again to the left
partition. From there, due to the shape of the related eigenvectors, they are
mapped again to the right partition, and so forth. This iterative process
leads to a WQA, which is the $\omega$-limit set of the iterations of segment
$PP_{1},$ i.e., the $\omega$-limit set of $(T_{1})^{n}(PP_{1})$, for
$n\rightarrow\infty.$

Its structure may appear quite weird. In fact, since segment $PP_{1}$ belongs
to an invariant region, all of its iterates remain in that region (divergence
is not possible). Note that linear homogeneous functions map segments
belonging to straight lines through the origin into segments belonging to
straight lines through the origin (Lemma 1). Thus, applying $T_{1}$ to segment
$PP_{1}$, in a finite number of iterations segment $(T_{1})^{k}(PP_{1})$
crosses the discontinuity line, leading to $2$ segments. Each of these
segments is then iterated similarly, each one (after a different number of
iterations) is crossing the discontinuity line, and thus a further division
occurs, leading to $2^{2}$ segments. This process continues indefinitely; the
iterates of the original segment generate $2^{n}$ segments, for any $n$. No
point can be mapped into itself in a finite number of iterations since cycles
do not exist, and all the segments belonging to $(T_{1})^{n}(PP_{1})$ are on
straight lines through the origin. The $\omega$-limit set of 
$(T_{1})^{n}(PP_{1})$ for $n\rightarrow\infty$ gives the attractor, a WQA. As in the
example in Fig.7(b), this attractor belong to an invariant area, a polygon
bounded by a finite number of critical segments. These critical segments are
the images of the segment on the discontinuity line included in the area
(shown in azure in Fig.8(a)).

Similar examples of WQAs are obtained for points in the wide red region of the
parameter plane shown in Fig.7(a).

\subsection{Further 2D PWL maps}

Next, we examine examples of maps defined by Jacobian matrices that differ
from those given in (\ref{MapT}), while keeping the same partitions, with
discontinuity set $x=-1$. Let us consider a triangular map with two different
linear functions having a common eigenvector, as in the following map:
\begin{equation}
T_{2}=\left\{
\begin{array}
[c]{c}%
T_{L}:X^{\prime}=J_{L}X \qquad \text{for } x<h, \quad J_{L}=\left[
\begin{array}
[c]{cc}%
a_{l} & b_{l}\\
0 & d_{l}%
\end{array}
\right] \\ \ \\ 
T_{R}:X^{\prime}=J_{R}X \qquad \text{for } x>h, \quad J_{R}=\left[
\begin{array}
[c]{cc}%
a_{r} & b_{r}\\
0 & d_{r}%
\end{array}
\right]
\end{array}
\right.  \ \ h=-1 \label{Jreal}%
\end{equation}
Both linear functions have eigenvector $y=0$ associated with eigenvalues
$a_{l}$ and $a_{r}$. The restriction of map $T_{2}$ on $y=0$ leads to the 1D
PWL circle map: 
$x^{\prime}=a_{l}x\ \ for\ x<-1$ and $x^{\prime}=a_{r}x\ \ for\ x>-1.$ 
That is, on eigenvector $y=0,$ the dynamics are those
occurring in map $F$ in (\ref{MapF}) with discontinuity point $x=-1$. Clearly,
the global behavior in the phase plane depends on the other eigenvalues of the
two linear functions, and their associated eigenvectors. Two examples are
presented in Figs.9,10.

\begin{figure}
	\begin{center}
		\includegraphics[width=0.8\textwidth]
		{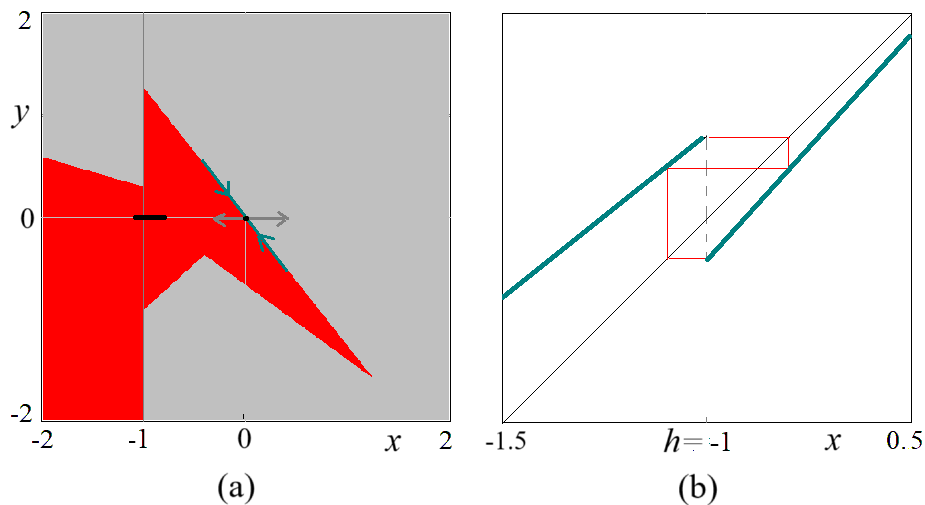} 
		\caption{\label{f9} \small{In (a), phase plane of map $T_{2}$ in (\ref{Jreal}) at $a_{l}=0.8,$
			$b_{l}=2,$ $d_{l}=0.9,$ $a_{r}=1.1,$ $b_{r}=1.5$ and $d_{r}=-0.8.$ The origin
			is a saddle fixed point, and the only attractor is the segment on eigenvector
			$y=0.$ The first return map on that segment is shown in Fig.9(b), a PWL circle
			map.}}
	\end{center}
\end{figure}
\begin{figure}
	\begin{center}
		\includegraphics[width=0.8\textwidth]
		{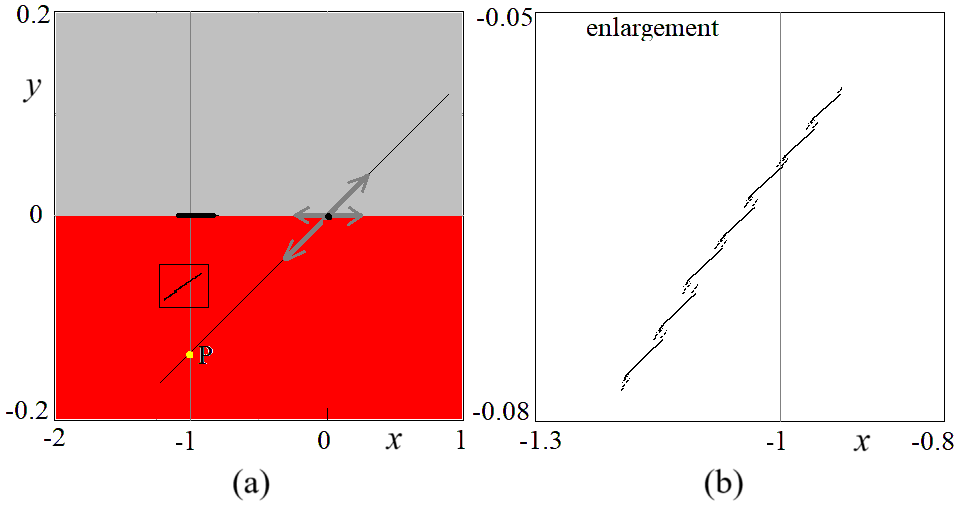} 
		\caption{\label{f10} \small{In (a), phase plane of map $T_{2}$ in (\ref{Jreal}) with the same
			parameter values as in Fig.9, except for $d_{r}=1.3.$ The origin is now a
			repelling node, and the only attractor is a WQA. An enlargement of the
			attractor is shown in Fig.10(b). The invariant segment on $y=0$ belongs to the
			border of the basin of attraction of the WQA.}}
	\end{center}
\end{figure}

In Fig.9, we consider the case where the fixed point $O$ is a virtual
attracting node, with both eigenvalues positive. In the right partition, the
real fixed point $O$ is a saddle with eigenvalue $\lambda_{1}=a_{r}=1.1,$
related to the eigenvector $y=0,$ and eigenvalue $\lambda_{2}=d_{r}=-0.8,$
related to the eigenvector with slope $s_{2}\simeq-1.267$. Figure 9(a) shows
that the only attractor is a segment on eigenvector $y=0,$ and the first
return map on that segment is a PWL circle map 
$x^{\prime}=a_{l}x=0.8x\ \ for\ x<-1$ and $x^{\prime}=a_{r}x=1.1x\ \ for\ x>-1,$ shown in
Fig.9(b). However, the basin of attraction depends on the global dynamics.
Since the fixed point $O$ is a saddle, we may expect divergent dynamics in
some regions of the phase plane. The gray region in Fig.9(a) represents
trajectories that eventually diverge. The boundary between the two basins (the
basin of divergent trajectories and the basin of the unique
attractor) consists of segments of the stable set of the origin (the
eigenvector associated with $\lambda_{2}=-0.8$), and segments of the
discontinuity line and their preimages.

In Fig.10, we modify only parameter $d_{r}$ (the second eigenvalue in the
right partition) so that the fixed point $O$ becomes a repelling node, with
positive eigenvalues. The restriction of the map to the eigenvector $y=0$
remains unchanged. The invariant interval on which the restriction of the map
is a PWL circle map remains the same as in the previous example (it is as in
Fig.9(b)). However, the global dynamics now differs. Since the fixed point
$O$ is a repelling node, eigenvector $y=0$ is repelling. Moreover, the
unstable eigenvector associated with eigenvalue $d_{r}=1.3$, with slope
$s_{2}\simeq0.1333,$ has the upper branch going to infinity, while the
opposite branch intersects the discontinuity line at a point $P$, so that a
segment $P_{-1}P$ on that eigenvector in the right partition is mapped into
$PP_{1}$ in the left partition. The iterates of this segment converge to a
WQA. This attractor is shown in Fig.10(a), highlighted by a rectangle, whose
enlargement is shown in Fig.10(b). The invariant segment on eigenvector $y=0$,
whose first return map is shown in Fig.9(b), lies on the boundary of the
basin of attraction (in red). Gray points denote divergent trajectories. 

Differently, if we consider the map with triangular functions defined as
follows:
\begin{equation}
T_{3}=\left\{
\begin{array}
[c]{c}%
T_{L}:X^{\prime}=J_{L}X \qquad \text{for } x<h, \quad  J_{L}=\left[
\begin{array}
[c]{cc}%
a_{1} & b_{1}\\
0 & c_{1}%
\end{array}
\right] \\
T_{R}:X^{\prime}=J_{R}X \qquad \text{for } x>h, \quad J_{R}=\left[
\begin{array}
[c]{cc}%
a_{2} & 0\\
b_{2} & c_{2}%
\end{array}
\right]
\end{array}
\right.  \ \ h=-1 \label{triang-2}%
\end{equation}
we have that $J_{L}$ has the eigenvector $y=0$ associated with eigenvalue
$a_{1}$, while $J_{R}$ has the eigenvector $x=0$ associated with eigenvalue
$c_{2}$, so that the functions in the two partitions do not have a common
eigenvector. For this discontinuous map, with triangular functions, we expect
the existence of WQAs. Two examples are shown in Fig.11(a,b), at parameter
values where the fixed point $O$ is a saddle for the function in the right
partition.

\begin{figure}
	\begin{center}
		\includegraphics[width=1\textwidth]
		{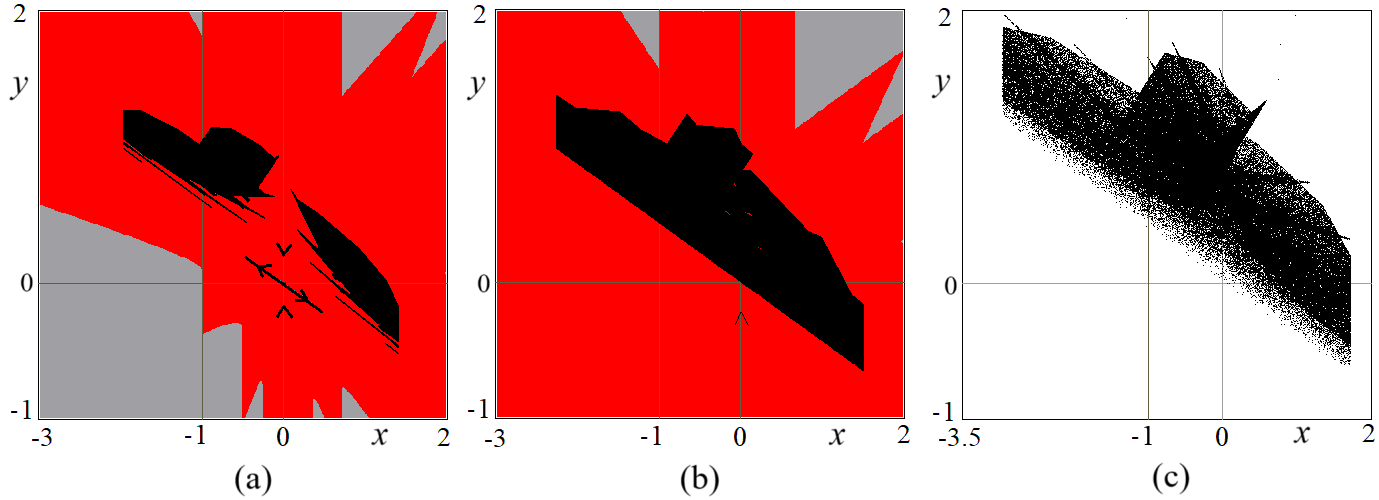} 
		\caption{\label{f11} \small{Phase plane of map $T_{3}$ in (\ref{triang-2}). In (a), at $a_{1}=1.1,$
			$b_{1}=1,$ $c_{1}=0.9,$ $a_{2}=-1.4,$ $b_{2}=1$ and $c_{2}=0.8.$ In (b),
			$a_{1}=0.9$ and $a_{2}=-1.5,$ with all other parameters as in (a). In (c),
			$a_{1}=0.9$ and $a_{2}=-1.72,$ with all other parameters as in (a). The
			transient of a divergent trajectory remains as a ghost of the former attractor
			for many iterations before diverging.}}
	\end{center}
\end{figure}
\begin{figure}
	\begin{center}
		\includegraphics[width=0.8\textwidth]
		{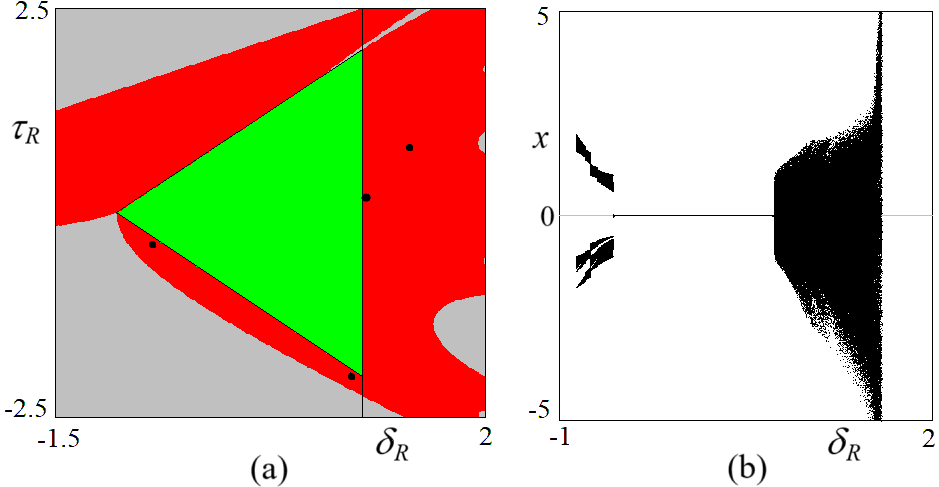} 
		\caption{\label{f12} \small{ In (a), 2D bifurcation diagram in the $(\delta_{R},\tau_{R})$ parameter
			plane for map $T_{4}$ in (\ref{NFprop}), with $\alpha=0.5.$ In (b), 1D
			bifurcation diagram as a function of $\delta_{R}$ for fixed $\tau_{R}=-0.5,$
			and $\alpha=0.5.$}}
	\end{center}
\end{figure}

In the case shown in Fig.11(a), the origin is a saddle (virtual) for both
partitions. In Fig.11(b,c), the origin is a virtual attracting node for
$T_{L}.$

The existence of a WQA can be explained as in previous cases. In the right
partition, the saddle fixed point has the stable eigenvector given by the
vertical line $x=0$. The branch of the unstable set of the fixed point $O,$
reaching the discontinuity line at a point $P,$ is such that $T_{R}$ maps a
segment, say $P_{-1}P,$ on that eigenvector into a segment $PP_{1}$ in the
$L$\ partition, where map $T_{L}$ applies. The WQA is the $\omega$-limit set
of $(T_{3})^{n}(PP_{1})$, for $n\rightarrow\infty.$

In both the cases shown in Fig.11(a,b), the basin of attraction of the WQA is
bounded by segments of the discontinuity line and segments of the eigenvector
$x=0$, along with their preimages. Gray points denote divergent trajectories.
The disappearance of the WQA occurs when there is a contact between the WQA
and its basin boundary.

It is worth noting that after the contact, all trajectories are divergent,
except for the fixed point $O$ (no other repelling cycle can exist). However,
near the bifurcation point, the system retains a form of "memory" of the
pre-existing WQA. At parameters just beyond the contact bifurcation, for
initial conditions close to the fixed point $O$, a long transient can be
observed as a "ghost" of the former WQA, before the trajectory eventually
diverges. An example of this behavior is illustrated in Fig.11(c). 

Considering the map in (\ref{MapT}) defined in the same regions but with the
following Jacobian matrices:%
\begin{equation}
T_{4}=\left\{
\begin{array}
[c]{c}%
T_{L}:X^{\prime}=J_{L}X \qquad \text{for } x<h, \quad  J_{L}=\left[
\begin{array}
[c]{cc}%
\alpha\tau_{R} & 1\\
-\alpha^{2}\delta_{R} & 0
\end{array}
\right] \\ \ \\ 
T_{R}:X^{\prime}=J_{R}X \quad \text{for } x>h, \quad  J_{R}=\left[
\begin{array}
[c]{cc}%
\tau_{R} & 1\\
-\delta_{R} & 0 
\end{array}
\right]
\end{array}
\right.  \ \ h=-1 \label{NFprop}%
\end{equation}
the linear maps in the two partitions have \textit{proportional eigenvalues,}
\textit{but not the same eigenvectors.} This property allows for the existence
of WQAs, as shown in the following examples. These attractors may occur both
when the eigenvalues of the unstable fixed point $O$ are real and when they
are complex conjugate.

Figure 12(a) shows the 2D bifurcation diagram in the parameter plane
$(\delta_{R},\tau_{R})$ for map $T_{4}$ in (\ref{NFprop}), at $\alpha=0.5.$
The stability triangle of the real fixed point $O$ is well evidenced, while
the red region denotes the existence of WQAs. The 1D bifurcation diagram
reported in Fig.12(b), as a function of $\delta_{R}$ at fixed 
$\tau_{R}=-0.5,$ clearly evidences the existence of WQAs, occurring both for
$\delta_{R}<1$ and for $\delta_{R}>1.$ The dynamics in the phase plane at the
four black dots marked in Fig.12(a) are shown in Fig.13, for $\delta_{R}<1,$
with real eigenvalues in the right partition, and in Fig.14, for 
$\delta_{R}>1,$ where the fixed point $O$ is a repelling focus.

When the eigenvalues are real, the mechanism leading to the existence of a WQA
is the same as described in the previous examples. In both cases in Fig.13,
the fixed point $O$ is a saddle, while it is a virtual attracting node for the
left partition. As a result, a segment related to one eigenvector of $O$ is
responsible for the crossing of the discontinuity line, entering the left
partition. Since the function in the left partition has a virtual attracting
fixed point, the dynamics converge to a WQA.

\begin{figure}
	\begin{center}
		\includegraphics[width=0.8\textwidth]
		{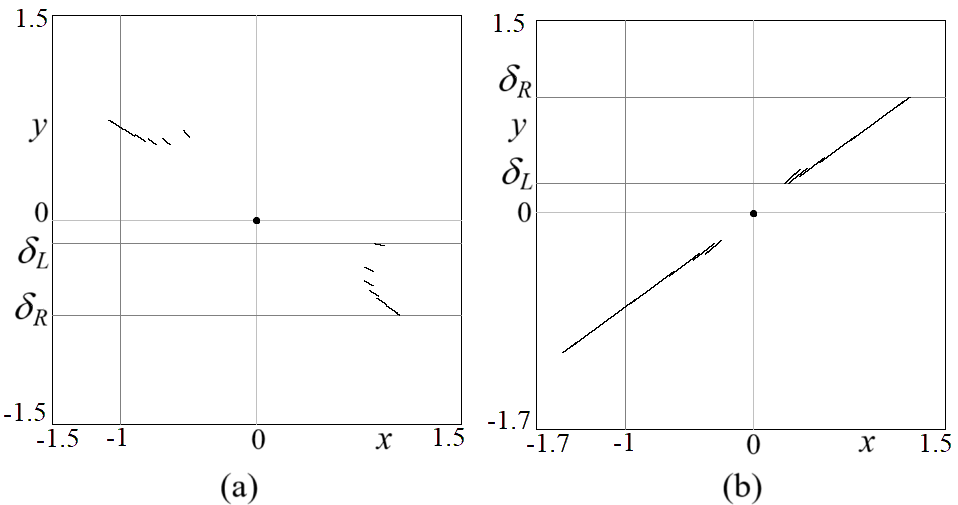} 
		\caption{\label{f13} \small{Map $T_{4}$ in (\ref{NFprop}), with $\alpha=0.5.$ In (a), WQA at
			$(\delta_{R},\tau_{R})=(-0.7,-0.37).$ The fixed point $O$ is a saddle, with
			eigenvalues of opposite sign, for the functions in both partitions. In (b),
			WQA at $(\delta_{R},\tau_{R})=(0.9,-1.96).$ The fixed point $O$ is a saddle,
			with two negative eigenvalues, for the functions in both partitions.}}
	\end{center}
\end{figure} 
\begin{figure}
	\begin{center}
		\includegraphics[width=0.8\textwidth]
		{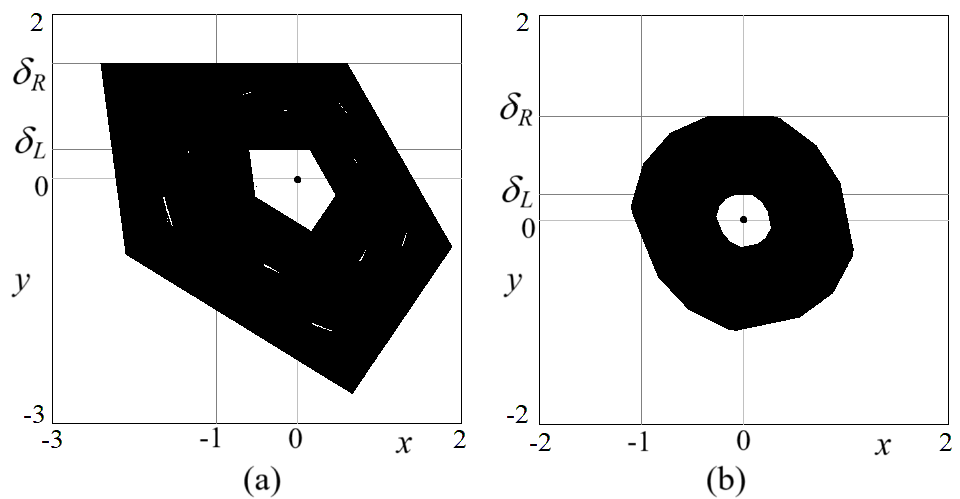} 
		\caption{\label{f14} \small{Map $T_{4}$ in (\ref{NFprop}) with $\alpha=0.5.$ In (a), WQA at
			$(\delta_{R},\tau_{R})=(1.4,0.8).$ In (b), WQA at 
			$(\delta_{R},\tau_{R})=(1.01,-0.2).$}}
	\end{center}
\end{figure}

In the case of complex eigenvalues, the mechanism leading to a WQA differs
from the previous one. In both examples shown in Fig.14, the fixed point $O$
is a repelling focus for the $R$ partition, while it is a virtual attracting
focus for the $L$ partition. Since divergence cannot occur, points from the
left partition are mapped to the right, from which they are repelled to the
left, and so forth. This mechanism leads to the WQA. In both examples of
Fig.14, the WQAs belong to an invariant area that can be determined by a
finite number of images of a segment of the discontinuity line (the segment
that is crossed by the attractor, similarly to the cases shown in Fig.7(b) and
Fig.8(a)).

\section{2D PWL homogeneous maps with different discontinuity sets}

\subsection{Discontinuity sets I: straight lines}

We first consider a map with the discontinuity line $y=x+h$, where $h=1.$ We
assign index $R$ to the region below the line, and $L$ to the region above it:%
\begin{equation}
T_{5}=\left\{
\begin{array}
[c]{c}%
T_{L}:X^{\prime}=J_{L}X \qquad  \text{for } y>x+h, \quad J_{L}=\left[
\begin{array}
[c]{cc}%
\tau_{L} & 1\\
-\delta_{L} & 0
\end{array}
\right] \\ \ \\ 
T_{R}:X^{\prime}=J_{R}X \qquad  \text{for } y<x+h, \quad J_{R}=\left[
\begin{array}
[c]{cc}%
\tau_{R} & 1\\
-\delta_{R} & 0
\end{array}
\right]
\end{array}
\right.  
\label{Reg-2}%
\end{equation}

In Fig.15(a), we show an example of a 2D bifurcation diagram in the parameter
plane $(\tau_{L},\tau_{R})$ for map $T_{5}$ in (\ref{Reg-2}), at fixed
$\delta_{R}=0.9$ and $\delta_{L}=0.8.$ Two cases with WQAs are shown in
Fig.15(b,c), corresponding to the black dots marked in Fig.15(a). 

\begin{figure}[h]
	\begin{center}
		\includegraphics[width=1\textwidth]
		{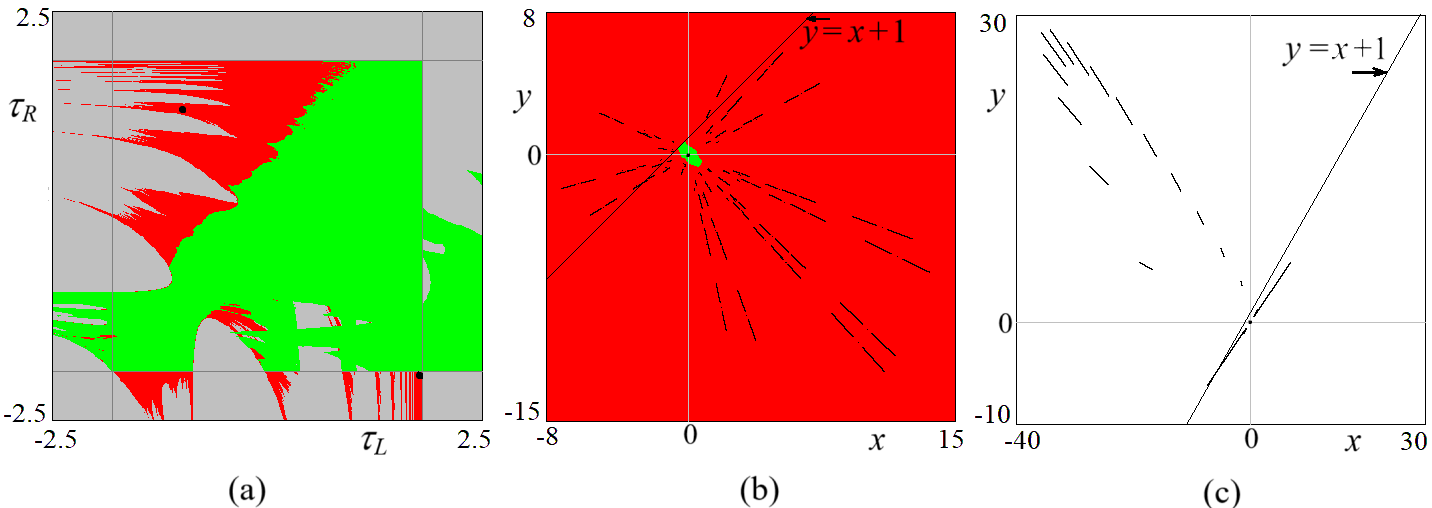} 
		\caption{\label{f15} \small{In (a), 2D bifurcation diagram in the $(\tau_{L},\tau_{R})$ parameter
			plane for map $T_{5}$ in (\ref{Reg-2}), with fixed $\delta_{R}=0.9$ and
			$\delta_{L}=0.8.$ In (b), phase plane at the black dot in (a) for 
			$(\tau_{L},\tau_{R})=(-1,1.3)$. The attracting fixed point $O$ coexists with a WQA.
			In (c), phase plane at the black point in (a) for 
			$(\tau_{L},\tau_{R})=(1.75,-1.91)$. The fixed point $O$ is a saddle for the functions in
			both partitions. The unique attractor is a WQA.}}
	\end{center}
\end{figure}

In the example in Fig.15(b), the origin is an attracting focus for both linear
functions, and a WQA exists. Here, the boundary of the basin of attraction of
the fixed point $O$ (shown in green in Fig.15(b)) consists of a segment of the
discontinuity line and a finite number of its preimages. The basin of
attraction has no point in the $L$ partition, where the map has a virtual
attracting focus at the origin. Consequently, the points from the $L$
partition are mapped to the $R$ partition, and rotate back to the $L$
partition again, and so forth. Divergence does not occur; a bounded attracting
set must exist.

Differently, in the example in Fig.15(c), the origin is a saddle for both
linear functions. The appearance of the WQA results from the unstable
eigenvector of the origin entering the left partition, converging to a
WQA.

Another case illustrated with two examples is shown in Fig.16, where the map
has two discontinuity lines, $y=x+h$ and $y=x-h$, with $h=1.$ Now index $R$
refers to the partition between the two straight lines, and $L$ outside of
that. The map is defined as:
\begin{equation}
T_{6}=\left\{
\begin{array}
[c]{c}%
T_{L}:X^{\prime}=J_{L}X \qquad  \text{for } |y-x|>h, \quad  J_{L}=\left[
\begin{array}
[c]{cc}%
\tau_{L} & 1\\
-\delta_{L} & 0
\end{array}
\right] \\ \ \\ 
T_{R}:X^{\prime}=J_{R}X \qquad \text{for } |y-x|<h, \quad  J_{R}=\left[
\begin{array}
[c]{cc}%
\tau_{R} & 1\\
-\delta_{R} & 0
\end{array}
\right]
\end{array}
\right.  \label{Reg-doppia}%
\end{equation} 

\begin{figure}[h]
	\begin{center}
		\includegraphics[width=1\textwidth]
		{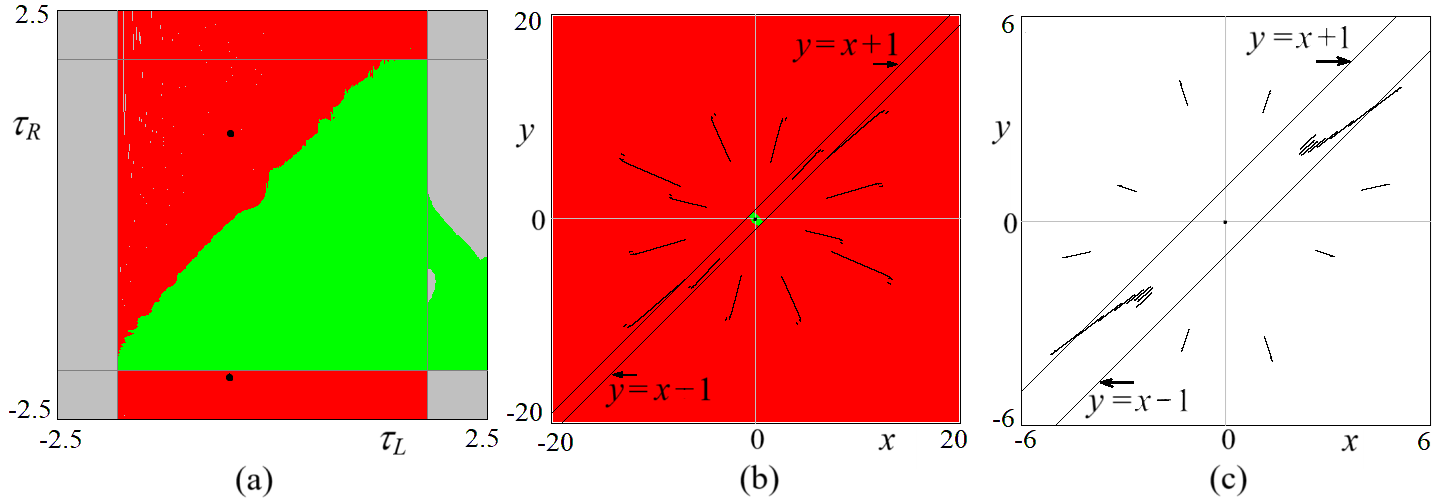} 
		\caption{\label{f16} \small{In (a), 2D bifurcation diagram in the $(\tau_{L},\tau_{R})$ parameter
			plane for map $T_{6}$ in (\ref{Reg-doppia}), with fixed $\delta_{R}=0.9$ and
			$\delta_{L}=0.8.$ In (b), phase plane at the black point in (a) for
			$(\tau_{L},\tau_{R})=(-0.5,1)$. The attracting fixed point $O$ coexists with
			a WQA. In (c), phase plane at the black dot in (a) for 
			$(\tau_{L},\tau_{R})=(-0.5,-1.95)$. The fixed point $O$ is a saddle for the $R$ partition
			and an attracting focus for the $L$ partition. The unique attractor is a WQA.}}
	\end{center}
\end{figure}

In Fig.16(a), we present an example of a 2D bifurcation diagram in the
parameter plane $(\tau_{L},\tau_{R})$ for map $T_{6}$ in (\ref{Reg-doppia}),
at fixed parameter values $\delta_{R}=0.9$ and $\delta_{L}=0.8.$ Two cases
with WQAs are shown in Fig.16(b,c), corresponding to the two black dots shown
in Fig.16(a). Map $T_{6}$ is now symmetric with respect to the fixed point
$O.$ It follows that an invariant set must be either symmetric with respect to
$O$, or the symmetric one also exists.

In Fig.16(b), the origin is an attracting focus for both functions $T_{R}$
and $T_{L}$ in (\ref{Reg-doppia}), while in Fig.16(c) the origin is a saddle
for $T_{R}$ and an attracting focus for $T_{L}.$ The unstable set of the fixed
point $O$ leads to a segment that enters the $L$ partition, converging to a WQA.

Similar results are observed when the two straight lines representing the
discontinuity sets are vertical. Examples of such cases are shown in
\cite{GRSSW-24}, \cite{GRSSW-25c}.

\subsection{Discontinuity sets II: circles}

Notably, it is not necessary to have straight lines as boundaries of the
partitions. For instance, let us consider a map where the discontinuity set is
defined by a circle, say $x^{2}+y^{2}=1.$ In this case, let us denote by index
$R$ the partition inside the circle, while $L$ corresponds to the region
outside. The map is given by:%
\begin{equation}
T_{7}=\left\{
\begin{array}
[c]{c}%
T_{L}:X^{\prime}=J_{L}X \qquad  \text{for } x^{2}+y^{2}>1, \quad  J_{L}=\left[
\begin{array}
[c]{cc}%
\tau_{L} & 1\\
-\delta_{L} & 0
\end{array}
\right] \\ \ \\ 
T_{R}:X^{\prime}=J_{R}X \qquad \text{for } x^{2}+y^{2}>1, \quad  J_{R}=\left[
\begin{array}
[c]{cc}%
\tau_{R} & 1\\
-\delta_{R} & 0
\end{array}
\right]
\end{array}
\right.  \label{Reg-3}%
\end{equation}

An example of a 2D bifurcation diagram for map $T_{7}$ in (\ref{Reg-3}), in
the parameter plane $(\delta_{L},\tau_{L}),$ is shown in Fig.17(a). Parameters
$(\delta_{R},\tau_{R})$ are fixed at $(1.4,0.8),$ so that the fixed point $O$
is a repelling focus. In this figure, the stability triangle of the virtual
fixed point for the left partition is colored in red, meaning that bounded
attractors exist. Two examples with WQAs are shown in Fig.17(b,c), at the
parameters represented by black dots in Fig.17(a). 

\begin{figure}[h]
	\begin{center}
		\includegraphics[width=1\textwidth]
		{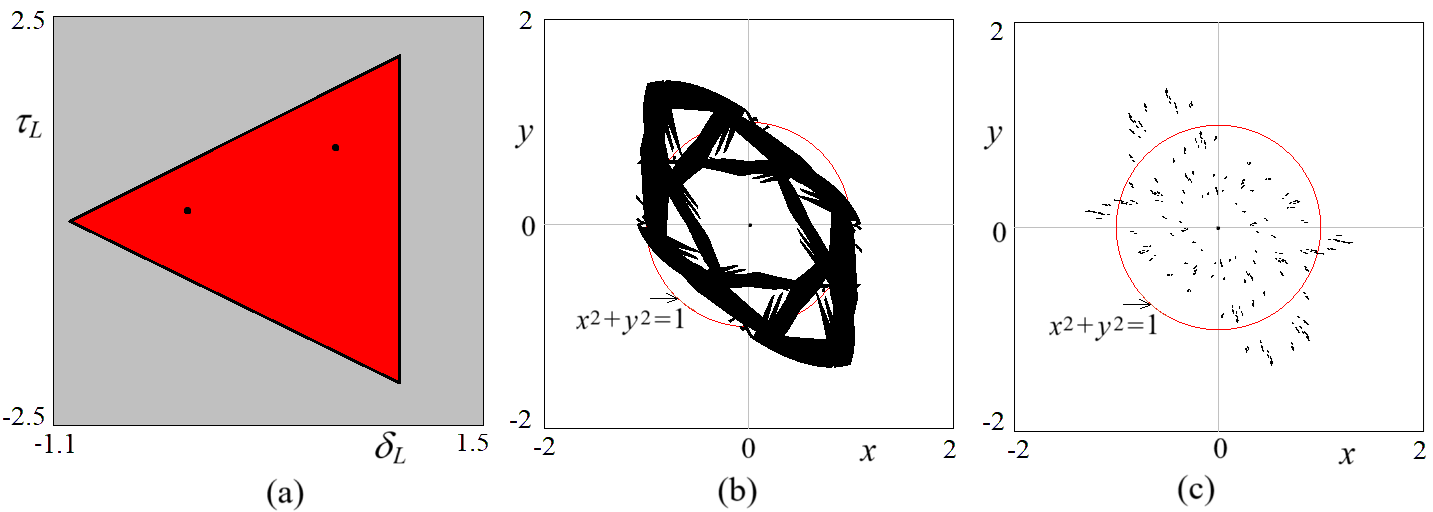} 
		\caption{\label{f17} \small{In (a), 2D bifurcation diagram in the $(\delta_{L},\tau_{L})$ parameter
			plane for map $T_{7}$ in (\ref{Reg-3}), with 
			$(\delta_{R},\tau_{R})=(1.4,0.8).$ WQAs are shown in (b) and (c). The discontinuity set is the red
			circle. In (b), phase plane at $(\delta_{L},\tau_{L})=(0.6,0.9),$ where the
			virtual fixed point $O$ is an attracting focus. In (c), phase plane at
			$(\delta_{L},\tau_{L})=(-0.3,0.1),$ where the virtual fixed point $O$ is an
			attracting node with eigenvalues of opposite signs.}}
	\end{center}
\end{figure}

In the first example, in Fig.17(b), map $T_{L}$ has complex eigenvalues. The
absence of divergent trajectories and the dynamics that rotate from inside the
circle (discontinuity set) to outside, and \textit{vice versa}, result in the
existence of a bounded attractor. In the second example, in Fig.17(c), map
$T_{L}$ has real eigenvalues; the parameters are inside the stability triangle
of $T_{L}$. The mechanism of formation of a WQA is similar to those observed
in previous examples. The eigenvector in the $L$ partition converging to the
virtual $O$ (inside the circle) leads to a segment in the $R$\ partition. The
iterates of this segment are forced to jump from inside the circle to outside,
and \textit{vice versa}, ultimately converging to the WQA. 

\textbf{Remark.} When the two matrices are proportional, for instance,
as described in Section 3.2 for map $T_{4}$ in (\ref{NFprop}), with $\alpha$
as proportional factor, but with different discontinuity sets, the map has
WQAs. In fact, the key property is that the two functions have proportional
eigenvalues but not the same eigenvectors. As a result, the generic attractor,
different from the fixed point $O$, is a WQA.

As an example, let us consider the same functions as in map $T_{4}$ in
(\ref{NFprop}), but with the unitary circle as discontinuity set, where
the $R$\ partition is inside the circle and $L$ outside. This gives the following map:
\begin{equation}
T_{8}=\left\{
\begin{array}
[c]{c}%
T_{L}:X^{\prime}=J_{L}X \qquad \text{for } x^{2}+y^{2}>1, \qquad  J_{L}=\left[
\begin{array}
[c]{cc}%
\alpha\tau_{R} & 1\\
-\alpha^{2}\delta_{R} & 0
\end{array}
\right] \\ \ \\ 
T_{R}:X^{\prime}=J_{R}X \qquad \text{for } x^{2}+y^{2}<1, \qquad J_{R}=\left[
\begin{array}
[c]{cc}%
\tau_{R} & 1\\
-\delta_{R} & 0
\end{array}
\right] \quad 
\end{array}
\right.  \label{T9}%
\end{equation}
Examples of WQAs for map $T_{8}$ in (\ref{T9}) are shown in Fig.18, with
$\alpha=0.5$. In Fig.18(a), the fixed point $O$ is a saddle (its eigenvectors
are shown inside the circle), while it is an attracting node for the $L$
partition. As in previous examples, the existence of a WQA may be related to
the eigenvectors. In Fig.18(b), the fixed point $O$ is a repelling focus,
while it is an attracting focus for the $L$ partition. As in other cases, the
existence of a WQA may be connected to the absence of divergent trajectories,
and the trajectories that from the $L$ partitions tend towards the virtual
origin, while from inside the circle, the trajectories are forced to rotate
going outside, eventually re-entering the $L$ partition, and so forth,
thereby sustaining a WQA.

\begin{figure}[h]
	\begin{center}
		\includegraphics[width=0.8\textwidth]
		{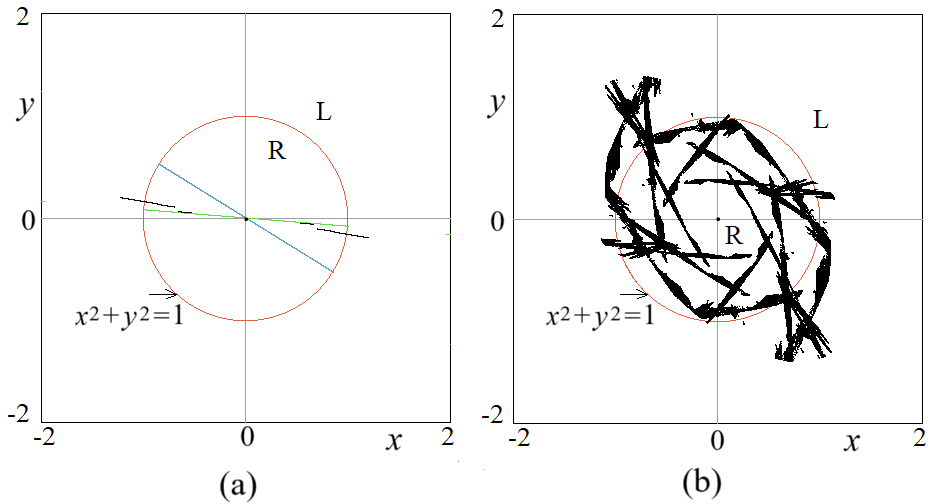} 
		\caption{\label{f18} \small{Phase plane with map $T_{8}$ in (\ref{T9}) with $\alpha=0.5.$ Two
			examples of WQAs. In (a), $\delta_{R}=0.2,$ $\tau_{R}=1.4.$ In (b),
			$\delta_{R}=1.4,$ $\tau_{R}=0.8.$}}
	\end{center}
\end{figure}

\subsection{2D PWL homogeneous map with the fixed point $O$ in two partitions}

In this subsection, we present several examples where the fixed point $O$ lies
on the border of two partitions, where the map is continuous, while an
additional discontinuity set is introduced.

We consider the previous examples, with a discontinuity set bounding the
partition, denoted as $D_{L}$ where the function $T_{L}$ is applied 
($D_{L}$ depends on the considered map). The remaining partition in the considered map,
denoted here as $D_{R},$ is then split in two regions, 
$D_{R}=R_{1}\cup R_{2},$ with $R_{1}$ for $x\geq0$ and $R_{2}$ for $x\leq0.$ The two Jacobian
matrices, in $R_{1}$\ and $R_{2}$, differ in terms of trace and determinant
following the standard 2D normal form structure. The map is defined as
follows:%
\begin{equation}
T_{9}=\left\{
\begin{array}
[c]{c}%
T_{L}:X^{\prime}=J_{L}X \qquad \text{for } X\in D_{L},  \qquad \qquad  J_{L}=\left[
\begin{array}
[c]{cc}%
\tau_{L} & 1\\
-\delta_{L} & 0
\end{array}
\right] \quad  \\
\\
T_{R2}:X^{\prime}=J_{R2}X \quad \text{for } X\in R_{2}(x\leq0), \quad  J_{R2}=\left[
\begin{array}
[c]{cc}%
\tau_{R2} & 1\\
-\delta_{R2} & 0
\end{array}
\right]  \\
\\
T_{R1}:X^{\prime}=J_{R1}X \quad \text{for } \in R_{1}(x\geq0), \qquad J_{R1}=\left[
\begin{array}
[c]{cc}%
\tau_{R1} & 1\\
-\delta_{R1} & 0
\end{array}
\right]
\end{array}
\right.  \label{T10}%
\end{equation}

In Fig.19, we present three examples of WQAs of map $T_{9}$ in
(\ref{T10}), each with a different discontinuity set. In all three examples,
the origin is a repelling focus for functions $T_{R1}$ and $T_{R2},$ while
it is an attracting focus for $T_{L}.$ Thus, the mechanism of formation of the
WQAs is always linked to trajectories spiraling outside the region that
includes the origin, while from outside the discontinuity set, from partition
$L$, the trajectories are spiraling towards the origin. The three cases differ
in the definition of $D_{L}.$ In Fig.19(a), $D_{L}$ corresponds to region
$x<-1,$ in Fig.19(b) to region $y>x+1,$ and in Fig.19(c) to region
$x^{2}+y^{2}>1$. 
\begin{figure}[h]
	\begin{center}
		\includegraphics[width=1\textwidth]
		{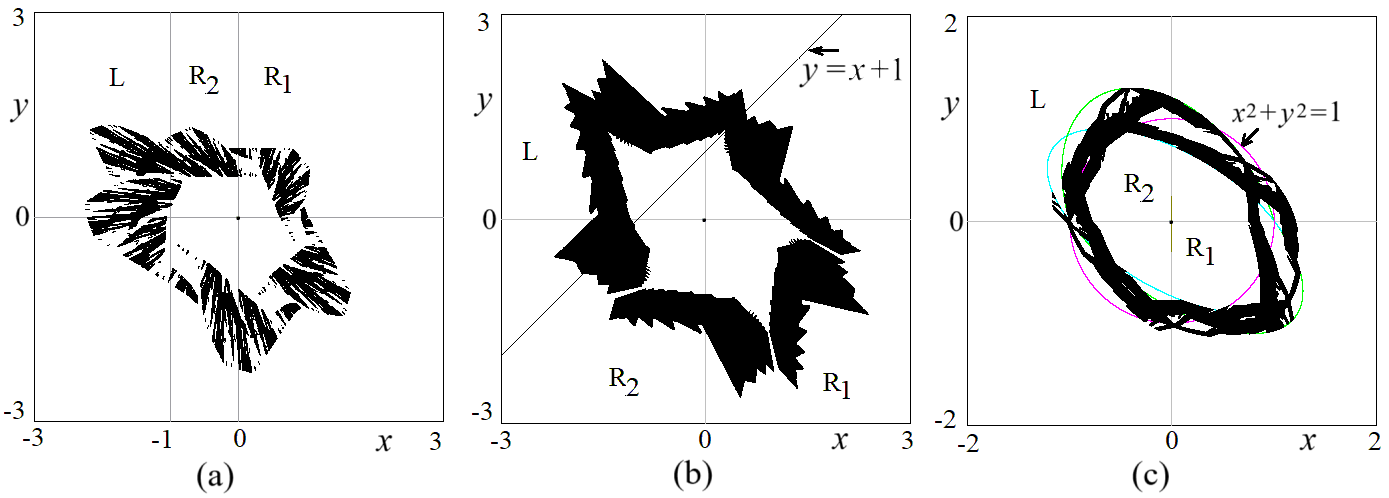} 
		\caption{\label{f19} \small{WQAs of map $T_{9}$ in (\ref{T10}) with three different partitions, and
				a single discontinuity set. In (a), the discontinuity set is line $x=-1$.
				Parameter values: $\tau_{R1}=0.8,$ 
				$\delta_{R1}=1.4,$ $\tau_{R2}=0.4,$ $\delta_{R2}=1.01,$ $\tau_{L}=0.9,$ $\delta_{L}=0.6.$ In (b), the
				discontinuity set is line $y=x+1$. Parameter values: 
				$\tau_{R1}=0.8,$ $\delta_{R1}=1.1,$ $\tau_{R2}=0.4,$ $\delta_{R2}=1.3,$ 
				$\tau_{L}=0.7,$ $\delta_{L}=0.9.$ In (c), the discontinuity set is circle 
				$x^{2}+y^{2}=1$. Parameter values: $\tau_{R1}=0.8,$ 
				$\delta_{R1}=1.1,$ $\tau_{R2}=0.4,$ $\delta_{R2}=1.3,$ $\tau_{L}=0.7,$ $\delta_{L}=0.9.$}}
	\end{center}
\end{figure}

Recall that the map is continuous and piecewise smooth in partition 
$R_{1}\cup R_{2},$ where it is the standard 2D normal form in the homogeneous case. Let
us denote this map as $T_{0},$ which, in our case, is applied at one side of
the discontinuity set. Map $T_{0}$ in the phase plane has been studied by
other authors, particularly in the conservative case, with 
$\delta_{R1}=\delta_{R2}=1$, see, e.g., 
\cite{Sivak,Beardon,Garcia,Lagarias-I,Lagarias-II,lagarias-III,Roberts}. In
that case, map $T_{0}$ depends on only two parameters, 
$(\tau_{R1},\tau_{R2}),$ and it was shown that its dynamics can be described by a circle map
with a well defined rotation number $\rho(\tau_{R1},\tau_{R2})$ that depends
on the parameters$.$ Hence, in the whole plane, the trajectories are either
all periodic with the same period (when $\rho$ is rational) or all
quasiperiodic and dense in closed curves that densely fill the plane (when
$\rho$ is irrational).

Thus, considering map $T_{9}$ with $\delta_{R1}=\delta_{R2}=1,$ for parameters
within suitable regions we may expect that these results hold for some points
in partition $R_{1}\cup R_{2}$. Let us show two examples with map $T_{0}$ in
partition $R_{1}\cup R_{2}$ of map $T_{9}$ in (\ref{T10}), considering the
unitary circle as discontinuity set, and region $D_{L}$(the $L$ partition)
outside the disc.

In Fig.20, the parameters of $T_{0}$ correspond to a rational rotation number
leading to 5-cycles and a region filled with 5-cycles coexists with a WQA.

\begin{figure}
	\begin{center}
		\includegraphics[width=0.8\textwidth]
		{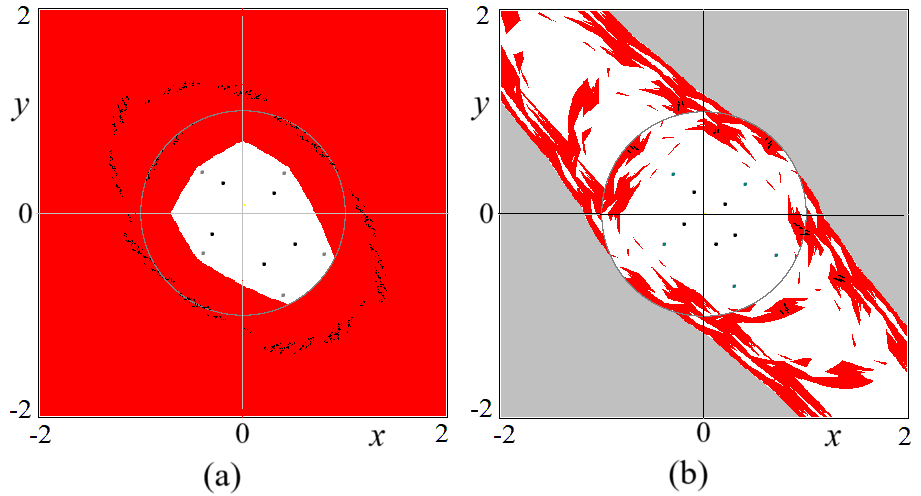} 
		\caption{\label{f20} \small{Phase plane of map $T_{9}$ in (\ref{T10}) and discontinuity set the
			unit circle. In (a), $\tau_{R1}=1,$ $\delta_{R1}=1,$
			$\tau_{R2}=0,$ $\delta_{R2}=1,$ $\tau_{L}=0.8,$ $\delta_{L}=0.98.$ The fixed
			point $O$ is a virtual attracting focus for map $T_{L}$ outside the circle.
			The white region is filled with 5-cycles, with two 5-cycles explicitly shown.
			Red region represents convergence to the WQA, shown by black
			points. In (b), map $T_{0}$ is the same as in (a), while for map $T_{L}$ the
			parameters are $\tau_{L}=0.3,$ $\delta_{L}=-0.8.$ The fixed point $O$ is now
			a virtual saddle for map $T_{L}$ outside the circle. Besides the conservative
			region filled with 5-cycles, the other points of white region are pre-periodic to a
			5-cycle. Gray region denotes divergence.}}
	\end{center}
\end{figure}

\begin{figure}
	\begin{center}
		\includegraphics[width=0.8\textwidth]
		{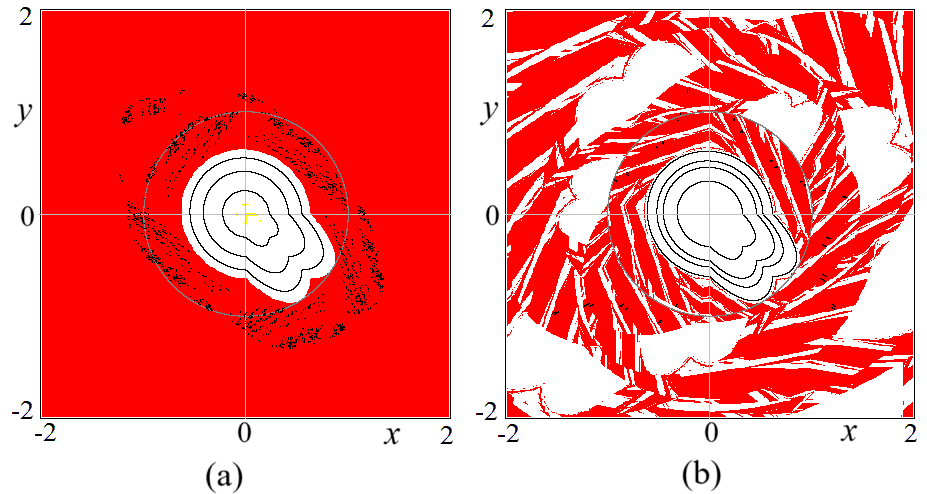} 
		\caption{\label{f21}\small{ Phase plane of map $T_{9}$ in (\ref{T10}), the unit circle is
			discontinuity set. In (a), $\tau_{R1}=1,$ $\delta_{R1}=1,$
			$\tau_{R2}=0.1,$ $\delta_{R2}=1,$ $\tau_{L}=0.8,$ $\delta_{L}=0.98.$ The
			fixed point $O$ is a virtual attracting focus for map $T_{L}$ outside the
			circle. The white region is filled with closed invariant curves, on which
			there are dense quasiperiodic trajectories. The red region represents
			convergence to the WQA, shown by black points. In (b), map $T_{0}$ remains the
			same as in (a), while for $T_{L}$ the parameters are 
			$\tau_{L}=-0.5,$ $\delta_{L}=0.8.$ 
			The fixed point $O$ is still a virtual attracting
			focus for $T_{L}$ outside the circle. Besides the conservative region
			filled with closed curves with quasiperiodic trajectories, the other points of white
			region are mapped to a closed curve in the conservative region.}}
	\end{center}
\end{figure}

Two different cases are shown. In Fig.20(a), map $T_{L}$ has complex
eigenvalues, with the origin a virtual attracting focus. The white region
(where map $T_{0}$ applies) represents the conservative region for map 
$T_{9}$. This region is bounded by an arc of discontinuity set and its images, and
it is filled with 5-cycles. The coexisting WQA has the basin in red. In
Fig.20(b), map $T_{L}$ has real eigenvalues (one positive, larger than $1$,
and one negative) and a region of divergent trajectories also exists. The
conservative region filled with 5-cycles is similar to the one in Fig.20(a).
However, in this case that region has also preimages outside. The union of all
the preimages results in the basin of attraction of the region filled with
5-cycles. That is, all the points in the white region of the figure, and
outside the conservative region, are pre-periodic to a 5-cycle. The basin of
the coexisting WQA is again shown in red.

In Fig.21, the parameters of $T_{0}$ correspond to an irrational rotation
number, leading to a suitable region filled with closed invariant curves on
which the trajectories are quasiperiodic, and there is coexistence with a WQA.
Two different cases are shown. In Fig.21(a), map $T_{L}$ has complex
eigenvalues, with the origin as a virtual attracting focus. The white region
(where map $T_{0}$ applies) is the conservative region for map $T_{9}$, filled
with closed invariant curves (with quasiperiodic trajectories). This region is
bounded by a closed invariant curve (the external one) that is tangent to the
discontinuity set, the unitary circle. The red region is the basin of
attraction of the WQA. Also in Fig.21(b), map $T_{L}$ has complex eigenvalues,
with the origin a virtual attracting focus, as in (a). However, now the
invariant conservative region filled with closed curves (that is similar to
the one in (a)) has other preimages. The white regions in the figure belong to
the basin of attraction of the invariant region. Their points are mapped into
one of the existing closed invariant curves. The red region belongs to the
basin of attraction of the coexisting WQA.

\section{Generalization to $\mathbb{R}^{n}$} 
The properties of 2D weird quasiperiodic attractors that we have described in
the previous sections can be generalized to the class of $n-$dimensional maps
for $n>2$. In this case, we call \textbf{nD weird quasiperiodic attractor} an
attractor $\mathcal{A}$ of an nD map $T$ that is a closed invariant set which
does not contain any periodic point (thus, it is neither an attracting cycle
nor a chaotic attractor). Moreover, the dynamics of $T$ on $\mathcal{A}$
cannot be studied by means of a first return map or by the restriction to a
set of lower dimension. In other words, an invariant set where the map is
reducible to a lower dimensional map, is not classified as nD WQA (although,
clearly, mD WQA, with $m<n$, are possible). 

We have already seen, in Lemma 1, the main properties of the maps in our
definition that hold in any dimension. We believe that these properties are
the essential elements needed to support the following

\begin{Conjecture} Let $T$ be an $nD$ map as given
in Definition 1. Then: 
\begin{itemize}[]

\item[(j)] A bounded $\omega$-limit set $\mathcal{A}$ 
different from the fixed point $O$ and from related local invariant
sets of $O$ when it is nonhyperbolic, can only be one of the following:
 
\begin{itemize}[]
\item[(ja)] a nonhyperbolic $k$-cycle, $k\geq2$ (this
occurs in $m$-dimensional sets, $m<n$, non intersecting any border and 
filled with cycles of the same symbolic sequence);

\item[(jb)] a finite number of $m$-dimensional sets, 
$m<n$, filled with quasiperiodic orbits; 

\item[(jc)] an invariant set, not structurally stable, 
on which the dynamics are reducible to a discontinuous $k$-dimensional map, $k<n$;

\item[(jd)] an $mD$ weird quasiperiodic attractor, where $2\leq m\leq n$.
\end{itemize} 

\item[(jj)] When no cycles exist filling $\mathcal{A}$ densely, then
 $\mathcal{A}$ exhibits (weak) sensitivity to initial 
conditions.
\end{itemize}
\end{Conjecture}

Let us consider a 3D example, with $X=(x,y,z).$ The simplest 3D map defined in
two partitions reads as follows:%
\begin{equation}
T=\left\{
\begin{array}
[c]{c}%
T_{L}:X^{\prime}=J_{L}X \qquad \text{for } x<h, \quad J_{L}=\left[
\begin{array}
[c]{ccc}%
\tau_{L} & 1 & 0\\
-\sigma_{L} & 0 & 1\\
\delta_{L} & 0 & 0
\end{array}
\right] \\
\\
T_{R}:X^{\prime}=J_{R}X \qquad \text{for } x>h, \quad  J_{R}=\left[
\begin{array}
[c]{ccc}%
\tau_{R} & 1 & 0\\
-\sigma_{R} & 0 & 1\\
\delta_{R} & 0 & 0
\end{array}
\right]
\end{array}
\right.  \label{T3D}%
\end{equation}
with $h\neq0$. Two examples of 3D WQAs are shown in Figs.22,23, for $h=-1$.

\begin{figure}
	\begin{center}
		\includegraphics[width=1\textwidth]
		{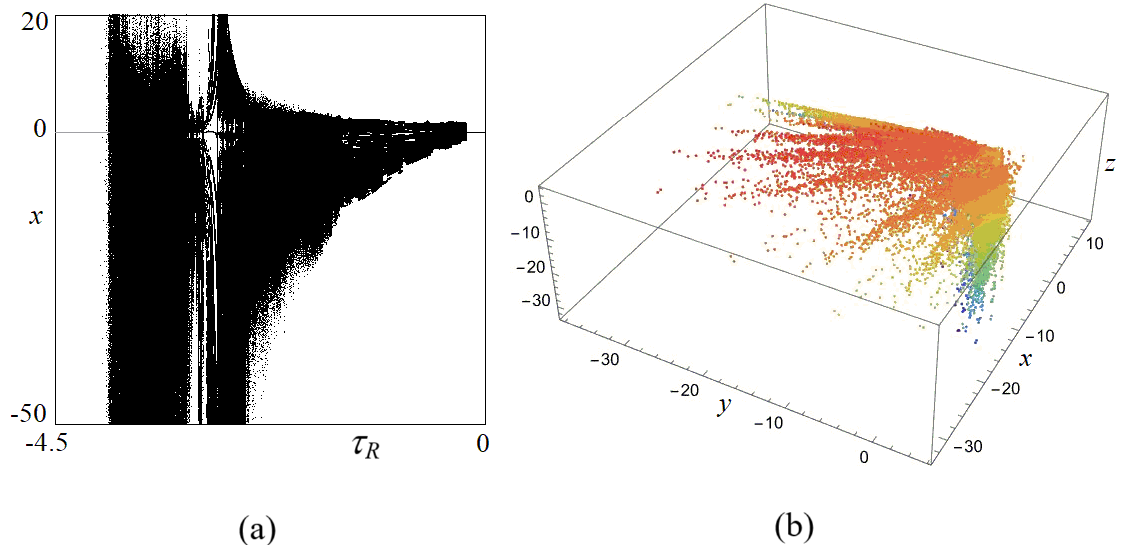} 
		\caption{\label{f22} \small{Map $T$ in (\ref{T3D}). In (a), 1D bifurcation diagram as a function of
			$\tau_{R}$ for $\sigma_{R}=0.2,$ $\delta_{R}=0.8,$ $\tau_{L}=0.3,$ 
			$\sigma_{L}=0.3,$ $\delta_{L}=0.9.$ The figure suggests the existence of WQAs. In
			(b), an example of a 3D WQA, at $\tau_{R}=-2.3$, with the other parameters as
			in (a).}}
	\end{center}
\end{figure}
\begin{figure}
	\begin{center}
		\includegraphics[width=1\textwidth]
		{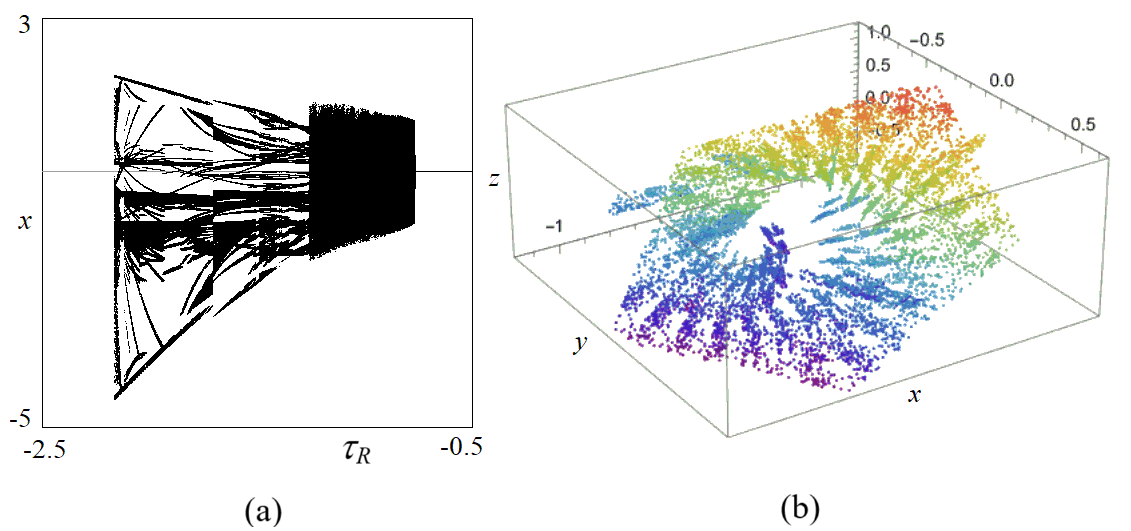} 
		\caption{\label{f23} \small{Map $T$ in (\ref{T3D}). In (a), 1D bifurcation diagram as a function of
			$\tau_{R}$ for $\sigma_{R}=-0.5,$ $\delta_{R}=0.9,$ $\tau_{L}=0.1,$ 
			$\sigma_{L}=-0.8,$ $\delta_{L}=0.5.$ In (b), an example of a 3D WQA, at 
			$\tau_{R}=-0.8$, with the other parameters as in (a).}}
	\end{center}
\end{figure}

In Fig.22(a), we present the 1D bifurcation diagram as a function of 
$\tau_{R}.$ Figure 22(b) shows an example of a 3D WQA. And in Fig.23(a), we present the 1D bifurcation diagram 
as a function of $\tau_{R}$ using a different set of parameters compared to Fig.22. Fig.23(b) shows
an example of a 3D WQA. Once again, such dynamics may easily be confused with chaotic motion. 

\section{Further research}
The maps belonging to the class considered in this work exhibit several
properties that deserve further investigation. In this section, we outline
some potential research directions. \medskip 

\textbf{(a) Shape of WQAs and sensitivity to parameter perturbation.}

In the class of discontinuous 2D PWL homogeneous maps (as in Definition 1),
the only structurally stable attractor different from the hyperbolic fixed
point $O$ is a weird quasiperiodic attractor. While a WQA persists under
parameter perturbation, its shape and structure need a deeper investigation.
Furthermore, shape and structure seem sensitive to parameter perturbation. In
some cases, even a small perturbation in a parameter can result in a drastic
change in the shape (i.e., a very different structure) of the WQA.\medskip

\textbf{(b) Existence of WQAs in a broader class of maps}.

We have shown that WQAs are generic attractors in a class of 2D discontinuous
PWL maps. However, we do not rule out the possible existence of WQAs in larger
classes of discontinuous piecewise smooth maps. That is, an attractor
$\mathcal{A}$ such that it does not include\ any periodic point and is the
$\omega$-limit set of quasiperiodic trajectories, that cannot be described by
means of a lower dimensional map (e.g., via a first return map), may exist in
other families of maps.\medskip

\textbf{(c) Maximum Lyapunov exponent in discontinuous maps.}

In the class of 2D PWL maps considered here, chaos cannot occur, leading us to
expect one negative Lyapunov exponent and one zero. However, the numerical
computations of the maximum Lyapunov exponent in discontinuous maps is often
not sufficiently reliable. Moreover, the maps in our class have attractors
showing sensitive dependence on initial conditions, due to the discontinuity
set. A particular scenario arises when the WQA has a structure almost dense in
some area. In such cases, a numerical computation of the maximum Lyapunov
exponent may yield a small positive value, despite the absence of a chaotic
behavior. This raises a question: could a numerical algorithm be developed for
discontinuous maps that differentiates between chaos and a WQA?\medskip

\textbf{(d) Transition from regular to chaotic dynamics.}

It is known that in the class of 1D PWL discontinuous Lorenz maps, the case of
a circle map denotes the transition from regular dynamics (in a gap map, where
chaos cannot occur) to chaotic dynamics (in an overlapping map). It would be
interesting to explore how the class of maps considered here behaves under
parameter perturbation, particularly under changes leading to different fixed
points (in place of a unique one) of the functions in the partitions. Is it
possible that such a perturbation leads a WQA to become a chaotic attractor?
We conjecture that this is possible, because in PWL maps, infinitely many
cycles, including homoclinic cycles, may appear under a small perturbation.\medskip

Let us consider the simplest case of map $T_{1}$ in (\ref{MapT}), with the
vertical straight line as discontinuity set, $h=-1$. We examine the attractors
that are numerically obtained when one of the linear functions is modified
into an affine function. Specifically, we consider the case shown in Fig.1(a)
(where the fixed point $O$ is repelling in the $R$ partition), and the
parameter point $(\tau_{L},\tau_{R})=(0.7,0.6).$ Thus, origin $O$ is a virtual
attracting focus for function $T_{L}$. At these parameter values, map $T_{1}$
in (\ref{MapT}) has a WQA, illustrated in Fig.24(a).

Now we consider the perturbed version of the map as follows:%
\begin{equation}
T_{1a}=\left\{
\begin{array}
[c]{c}%
T_{L}:\left\{
\begin{array}
[c]{c}%
x^{\prime}=\tau_{L}x+y+\mu_{L}\\
y^{\prime}=-\delta_{L}x \qquad \quad 
\end{array}
\right.  \quad \text{for } x<h, \quad J_{L}=\left[
\begin{array}
[c]{cc}%
\tau_{L} & 1\\
-\delta_{L} & 0
\end{array}
\right] \\ \ \\ 
T_{R}:\left\{
\begin{array}
[c]{c}%
x^{\prime}=\tau_{R}x+y\\
y^{\prime}=-\delta_{R}x \quad 
\end{array}
\right.  \qquad \quad \text{for } x>h, \quad  J_{R}=\left[
\begin{array}
[c]{cc}%
\tau_{R} & 1\\
-\delta_{R} & 0
\end{array}
\right]
\end{array}
\right.  \label{Perturb1}%
\end{equation}
where we add a constant term, $\mu_{L},$ in the definition of $x^{\prime}$ in
the $L$ partition. The qualitative shape of the attracting set remains
unchanged for both $\mu_{L}>0$ and $\mu_{L}<0,$ when close to $0$, as shown
in Fig.24 (for $\mu_{L}>0)$ and Fig.25 (for $\mu_{L}<0)$. However, now the
attractors may be chaotic.

\begin{figure}
	\begin{center}
		\includegraphics[width=1\textwidth]
		{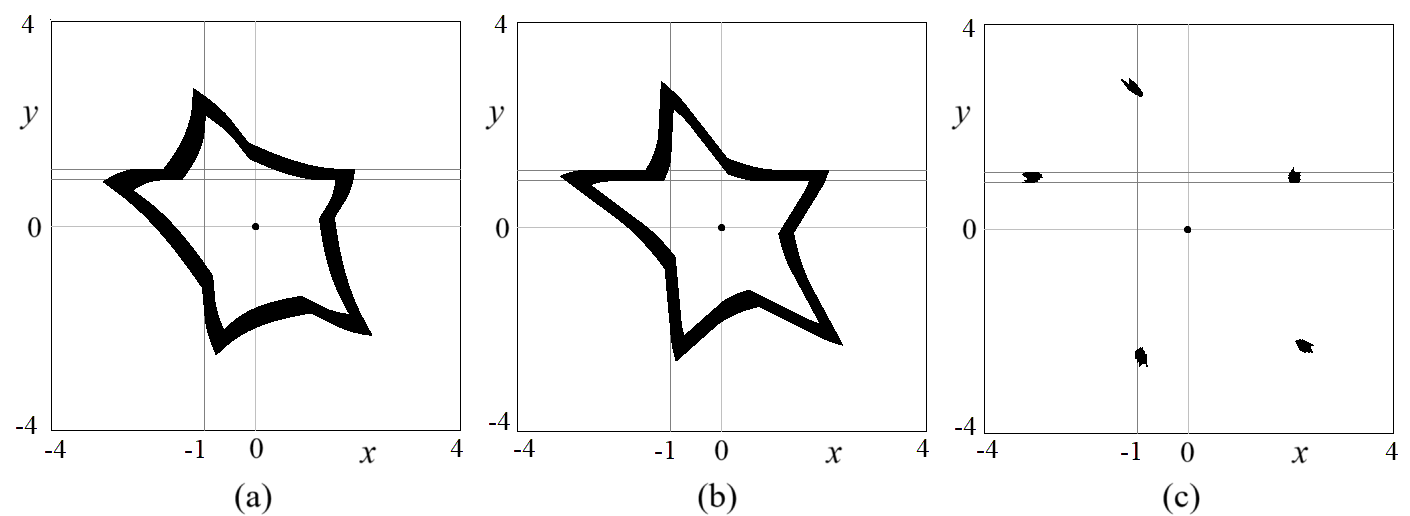} 
		\caption{\label{f24} \small{In (a), WQA of map $T_{1}$ in (\ref{MapT}), at $\tau_{L}=0.7,$
			$\delta_{L}=0.9,$ $\tau_{R}=0.6,$ $\delta_{R}=1.11.$ In (b) and (c), map
			$T_{1a}$ in (\ref{Perturb1}), with the same parameters as in (a) and 
			$\mu_{L}=0.03$ in (b), $\mu_{L}=0.05$ in (c).}}
	\end{center}
\end{figure}

\begin{figure}
	\begin{center}
		\includegraphics[width=1\textwidth]
		{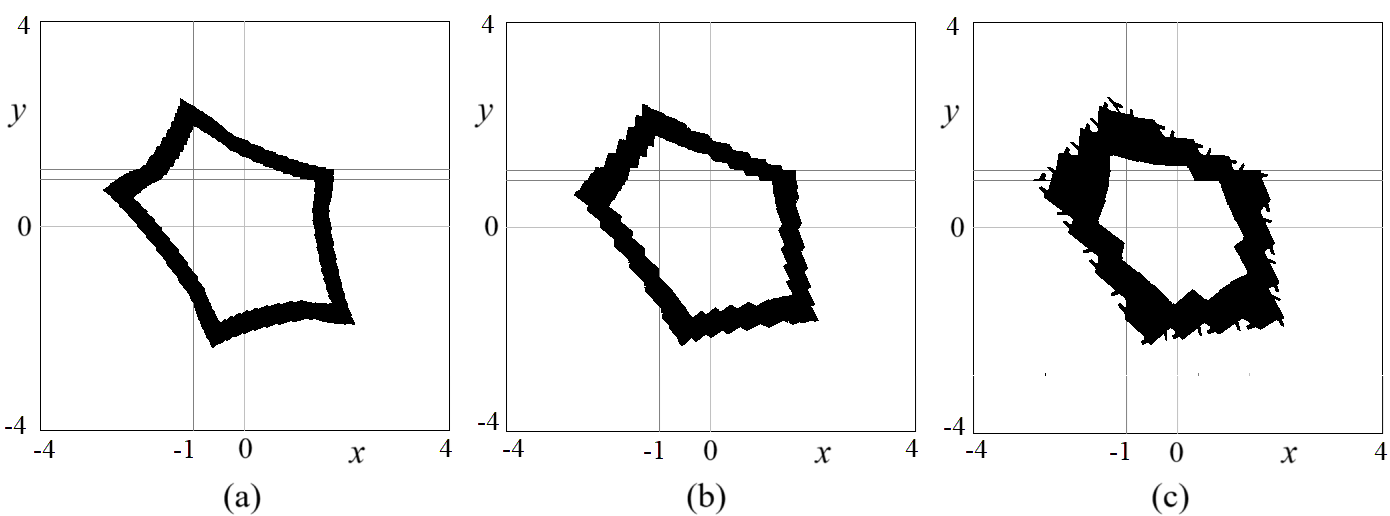} 
		\caption{\label{f25} \small{Map $T_{1a}$ in (\ref{Perturb1}), at $\tau_{L}=0.7,$ $\delta_{L}=0.9,$
			$\tau_{R}=0.6,$ $\delta_{R}=1.11.$ In (a), $\mu_{L}=-0.03.$ In (b), 
			$\mu_{L}=-0.1.$ In (c), $\mu_{L}=-0.3.$}}
	\end{center}
\end{figure}

Similar dynamic behaviors are obtained with map $T_{1}$ in (\ref{MapT}), when
a constant term $\mu_{R}$ is added to the function in the right partition, as
well as when constant terms are introduced in both functions.

\section{Conclusions} 
In this work, we have considered a class of $n-$dimensional piecewise linear
discontinuous maps, as defined in Definition 1, that generalizes the family of
maps considered in \cite{GRSSW-25a}. These maps can have\ a new kind of
attractor, called a weird quasiperiodic attractor (WQA). Piecewise smooth maps
and, in particular, PWL maps, are widely used in many applied fields. Our
initial interest in this class of maps arose from financial market models.

The characteristic property of the maps in our class, besides the
discontinuity, is that all the functions\ in the definition have the same real
fixed point. Although numerical simulations may sometimes suggest chaotic
behavior, we have shown that the dynamics associated with these attractors
cannot be chaotic.

The 1D case was previously considered in \cite{GRSSW-25b}. In this work, we
provided a detailed investigation of the 2D case, along with some generic
properties that hold for any dimension $n$ (Lemma 1). In particular, the
maps\ satisfying our definition cannot have hyperbolic cycles different from
the fixed point, which is sufficient to prove that a chaotic attractor cannot exist.
Cycles can only be nonhyperbolic, and are nongeneric, meaning they do not
persist under parameter perturbation.

In Section 2 we have shown that on straight lines through the fixed point invariant sets can exist, 
for 2D maps in our
definition, in which the restriction of the map can be reduced to a 1D map. In
such cases, Theorem 2 establishes that the map is necessarily related to a PWL
circle map, whose dynamics are well know, and depend on the rotation number,
rational or irrational. The main result is given in Theorem 3, where we show
that the dynamics on an invariant set are either reducible to those 
of a 1D map, or to a WQA. This conclusion holds independently of the (finite)
number and shape of the discontinuity sets, as well as of the location of the
real fixed point (which may be internal to a partition or on a border).
Moreover, a WQA may coexist with other attractors or invariant sets, as well
as with divergent trajectories.

In Section 3, we have illustrated our results through various examples, using
PWL maps with different Jacobian matrices but the same discontinuity set. In
Section 4, we have explored the effects of different discontinuity sets. A
generalization to nD WQA is possible, as noted in Section 5, where we have
presented numerical examples of 3D WQAs.

This new kind of attractor needs to be better investigated, and several
directions of future research have been outlined in Section 6. In particular,
the intrinsic structure of a WQA has yet to be understood, as well as the
mechanisms that may lead to its appearance/disappearance. While we have
provided some initial insights into these aspects, further exploration is
necessary to understand the properties of this new type of attractor.
\bigskip

\textbf{Acknowledgements}

Laura Gardini thanks the Czech Science Foundation (Project 22-28882S), the
VSB---Technical University of Ostrava (SGS Research Project SP2024/047), the
European Union (REFRESH Project-Research Excellence for Region Sustainability
and High-Tech Industries of the European Just Transition Fund, Grant
CZ.10.03.01/00/22 003/000004). Davide Radi thanks the Gruppo Nazionale di
Fisica Matematica GNFM-INdAM for financial support. The work of Davide Radi
and Iryna Sushko has been funded by the European Union - Next Generation EU,
Mission 4: "Education and Research" - Component 2: "From research to
business", through PRIN 2022 under the Italian Ministry of University and
Research (MUR). Project: 2022JRY7EF - Qnt4Green - Quantitative Approaches for
Green Bond Market: Risk Assessment, Agency Problems and Policy Incentives -
CUP: J53D23004700008.

\end{document}